%% file: main.tex
\documentclass[mathpazo]{cicp}

\usepackage{amsmath,amsfonts,amsthm,amssymb}
\usepackage{algorithm2e}
\usepackage{bm,dsfont}
\usepackage{booktabs}
\usepackage{graphicx}
\usepackage{array}
\usepackage{xcolor}
\usepackage{sidecap}
\usepackage{listings}
\usepackage{times}
\usepackage{url}
\usepackage{hyperref}

\ifpdf
  \DeclareGraphicsExtensions{.eps,.pdf,.png,.jpg}
\else
  \DeclareGraphicsExtensions{.eps}
\fi
\ifpdf
\hypersetup{
	pdftitle={On the convergence of PINNs},
	pdfauthor={Y. Shin, J. Darbon, G. E. Karniadakis}
}
\fi

\newtheorem{corollary}[theorem]{Corollary}
\newtheorem{claim}[theorem]{Claim}
\newtheorem{assumption}{Assumption}[section]

\newcommand{\x}{\textbf{x}}

\begin{document}
	
	\title{On the convergence of physics informed neural networks for linear second-order elliptic and parabolic type PDEs}

    \author[Shin Y et.~al.]{Yeonjong Shin\affil{1}, J\'er\^ome Darbon\affil{1}, and George Em Karniadakis\affil{1}}
    \address{
    \affilnum{1}\ Division of Applied Mathematics, Brown University, Providence, RI 02912, USA}
    \emails{{\tt yeonjong\_shin@brown.edu} (Y.~Shin), {\tt jerome\_darbon@brown.edu} (J.~Darbon),
        {\tt george\_karniadakis@brown.edu} (G.~E.~Karniadakis).}

	\begin{abstract}
		Physics informed neural networks (PINNs) are deep learning based techniques for solving partial differential equations (PDEs)
		encounted in computational science and engineering.
		Guided by data and physical laws, PINNs find a neural network that approximates the solution to a system of PDEs.
		Such a neural network is obtained by minimizing a loss function in which 
		any prior knowledge of PDEs and data are encoded.
		Despite its remarkable empirical success
		in one, two or three dimensional problems, there is little theoretical justification for PINNs.
		
		As the number of data grows, PINNs generate a sequence of minimizers which correspond to a sequence of neural networks.
		We want to answer the question: Does the sequence of minimizers converge to the solution to the PDE? 
		We consider two classes of PDEs: linear second-order elliptic and parabolic.
		By adapting the Schauder approach and the maximum principle, 
		we show that 
		the sequence of minimizers strongly converges to the PDE solution in $C^0$.
		Furthermore, we show that if each minimizer satisfies the initial/boundary conditions,
		the convergence mode becomes $H^1$.
		Computational examples are provided to illustrate our theoretical findings.
		To the best of our knowledge, this is the first theoretical work that shows the consistency of PINNs.
	\end{abstract}

	\ams{65M12, 41A46, 35J25, 35K20}
	
	\keywords{Physics Informed Neural Networks, Convergence, H\"{o}lder Regularization, Elliptic and Parabolic PDEs, Schauder approach}
	
	\maketitle
	
	\input introduction

	\input setup
	\input analysis

	\input examples

	\input conclusion
	
	\appendix
	\input appendix

	\input acknowledgement

	\bibliographystyle{plain}
	\bibliography{references}
\end{document}

%% file: introduction.tex
\section{Introduction}
Machine learning techniques using deep neural networks have been successfully applied in various fields \cite{Lecun_Nature15_DeepLearning}
such as computer vision and natural language processing.
A notable advantage of using neural networks 
is its efficient implementation
using a dedicated hardware
(see \cite{Lagaris_98_ANN-ODE-PDE, Darbon_19_NNsHJ}).
Such techniques have also been applied in 
solving partial differential equations
(PDEs) \cite{Raissi_19_PINNs, Lagaris_98_ANN-ODE-PDE, Dissanayake_94_ANN-PDE, Lagaris_00_ANN-Irregular, Berg_18_Unified, Sirignano_JCP18_DGM}, 
and it has become a new sub-field under the name of Scientific Machine
Learning (SciML) \cite{Baker_19_DCworkshop,Lu_19_Deepxde}.
The term Physics-Informed Neural Networks (PINNs) was introduced in \cite{Raissi_19_PINNs} and it has become one of the most popular deep learning methods in SciML. 
PINNs employ a neural network as a solution surrogate
and seek to find the best neural network guided by data and physical laws expressed as PDEs.

A series of works have shown the effectiveness of PINNs in one, two or three dimensional problems: fractional PDEs \cite{Pang_SISC19_fPINNs,Song_19_fPINNs}, stochastic differential equations \cite{Zhang_19_SPDE_PINNs,HanE_18_DLSPDEs}, 
biomedical problems \cite{Raissi_Nature20_HFM}, and fluid mechanics \cite{Mao_20_HighSpeedFlows}.
Despite such remarkable success in these and related areas, 
PINNs lack theoretical justification.
In this paper, we provide a mathematical justification of PINNs.

\begin{figure}[htbp]
	\centerline{
		\includegraphics[width=10cm]{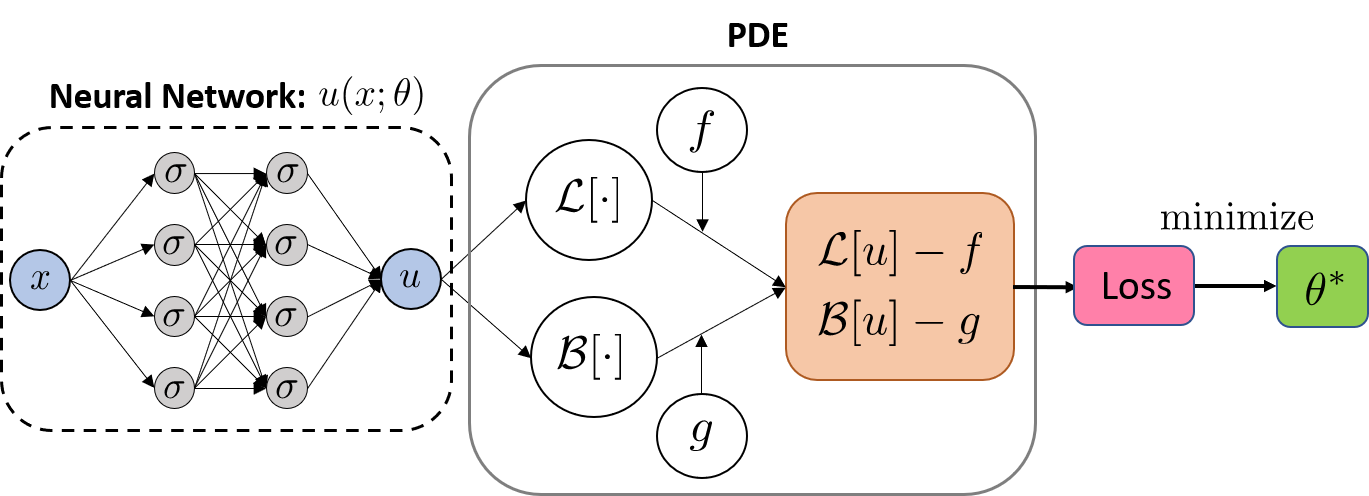}
	}
	\caption{
		 Schematic of physics-informed neural networks (PINNs).
		 The left part visualizes a standard neural network parameterized by $\theta$. 
		 The right part applies the given physical laws to the network.
		 $\mathcal{L}$ and $\mathcal{B}$ are the differential and the boundary operators, respectively.
		 The PDE data ($f, g$) are obtained from random sample points.
		 The loss function is computed by evaluating 
		 $\mathcal{L}[u]$ and $\mathcal{B}[u]$ on the sample points,
		 which can be done efficiently through automatic differentiation \cite{Baydin_17AD}.
		 Minimizing the loss with respect to the network's parameters $\theta$ produces 
		 a PINN $u(x;\theta^*)$, which serves as an approximation 
		 to the solution to the PDE.
	}
	\label{fig:PINNs_Schematic}
\end{figure}

One of the main goals of PINNs is to approximate the solution to the PDE.
For readers' convenience, we provide a concrete example to explain the PINNs method.
Let us consider a 1D Poisson equation on $U=(0,1)$ with the Dirichlet boundary conditions:
\begin{equation*}
    \mathcal{L}[u](x) = f(x), \quad \forall x \in U, \qquad \mathcal{B}[u](x) = g(x), \quad \forall x \in \partial U,
\end{equation*}
where the differential operator $\mathcal{L}$ is the Laplace $\Delta$ operator and the boundary operator $\mathcal{B}$ is the identity operator.
Suppose that 
the differential operator $\mathcal{L}$ and the boundary operator $\mathcal{B}$ 
are known to us, however, the PDE data (i.e., $f$ and $g$) are only known at some sample points.
In other words, although we do not know what $f$ and $g$ are exactly,
we can access them through their pointwise evaluations.
{Here and throughout the paper, we assume the high regularity setting for PDEs where point-wise evaluations are defined.}
We refer to a set of available pointwise values of $f$ and $g$ as the training data set.
The loss functionals
are then designed to penalize functions that fail to satisfy 
both governing equations (PDEs) and boundary conditions
\emph{on the training data}.
For example, if $(j/5, f(j/5))$ for $j=1,\dots,4$, and $(0,g(0)$ and $(1,g(1))$ are given to us as the training data, a prototype PINN loss functional \cite{Raissi_19_PINNs} is defined by
\begin{equation*}
    \text{Loss}(u) = \frac{1}{4}\sum_{j=1}^4 (\Delta u(j/4) - f(j/4))^2 + \frac{1}{2}\left[(u(0) - g(0))^2 + (u(1)-g(1))^2\right]. 
\end{equation*}
PINNs seek to find a neural network that minimizes the loss in a class of neural networks.
A minimizer serves as an approximation to the solution to the PDE.
{
Since neural networks
are parameterized with finitely many variables,}
the loss functional restricted to it becomes a function of network parameters, which is called the loss function.
With a slight abuse of terminology, we shall not distinguish the loss functional and the loss function and both will be simply called the loss function.
In Figure~\ref{fig:PINNs_Schematic}, we provide a schematic of PINNs.

PINNs are different approaches to the traditional variational principle that minimizes an energy functional \cite{Attouch_14Book_variational,Brezis_10Book_functional,Weinan_18DeepRitz}.
The most distinctive difference between them is that not all PDEs satisfy a variational principle, however, the formulation of PINNs does not require the considered PDE to have a variational principle
and can be applied to broad classes of PDEs as shown in aforementioned works.

As explained, given a set of $m$-training data, 
PINNs require one to choose  
a class of networks $\mathcal{H}_n$ and a loss function.
$n$ represents a complexity of the class, e.g., the number of parameters.
$\mathcal{H}_n$ may depend on the number of training data and also the data itself.
The goal is then to find a minimizer $h_m \in \mathcal{H}_n$ of the loss. 
The optimization process is referred to as training.
However, even in an extreme case where $\mathcal{H}_n$ contains the exact solution $u^*$ to PDEs
and a global minimizer is found, since there could be multiple (often infinitely many) global minimizers, there is no guarantee that a minimizer one found and the solution $u^*$ coincide.
By assuming an idealized setup of a global minimizer always being found, 
PINNs generate a sequence of neural networks with respect to the number of data.
We want to answer the question: Does the sequence of neural networks converge to the solution to PDEs?

\begin{figure}[htbp]
	\centerline{
		\includegraphics[width=12cm]{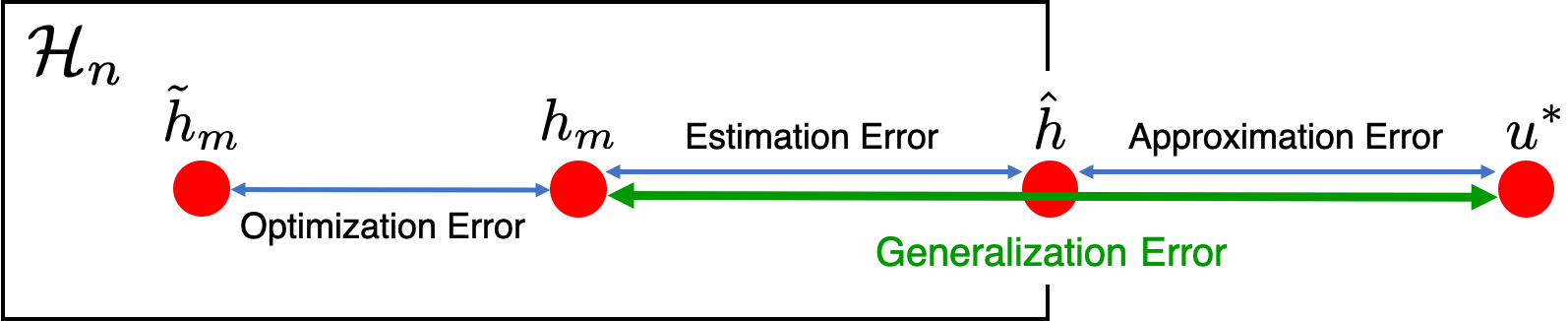}
	}
	\caption{
		Illustration of the total errors.
		$\mathcal{H}_n$ is the chosen function class.
		$u^*$ is the solution to the underlying PDE.
		The number of training data is $m$.
		$h_m$ is a minimizer of the loss with $m$ data.
		$\hat{h}$ is a function in $\mathcal{H}_n$
		that minimizes the loss with infinitely many data.
		$\tilde{h}_m$ is an approximation that one obtains in practice, e.g., 
		the result obtained after 1M epochs of a gradient-based optimization. 
	}
	\label{fig:Error}
\end{figure}

The total errors of neural networks-based supervised learning can be decomposed into three components \cite{niyogi1999generalization,bottou2008tradeoffs}: 
(a) approximation error, (b) optimization errors, and (c) estimation error.
We illustrate the decomposition of the total errors in Figure~\ref{fig:Error}.
%
The approximation error (a) is relatively well understood.
\cite{Pinkus_99_ATofMLP} showed that a single layer neural network with a sufficiently large width
can uniformly approximate a function and its partial derivative. 
It also has been shown that neural networks are capable of approximating
the solutions for some classes of PDEs: Quasilinear parabolic PDEs 
\cite{Sirignano_JCP18_DGM}, the Black-Scholes PDEs 
\cite{Grohs_18_BSPDEs}, and the Hamilton-Jacobi PDEs \cite{Darbon_19_NNsHJ, Darbon_20_NNsHJ}.
%
%
The optimization error (b) is, however, poorly understood as the objective function is highly nonconvex. Optimization often involves many engineering tricks and tedious trial and error type fine-tuning of parameters.
Gradient-based optimization methods are commonly used for the training.
Numerous variants of the stochastic gradient descent method
have been proposed \cite{Ruder_16_GDoverview}. 
In particular, \texttt{Adam} \cite{Kingma_14_Adam} and \texttt{L-BFGS} \cite{Liu_89_LBFGS} are popularly employed in the context of PINNs \cite{Lu_19_Deepxde}.
Empirically, the solutions found by gradient-based optimization are shown to 
perform well in various challenging tasks.
However, to the best of our knowledge, there is no guarantee that gradient-based optimization will find a global minimum for general machine learning problems including PINNs.
This is one of widely open problems.
The works of \cite{Jagtap_19LAAF,Wang_20_GradPathPINN} 
studied a couple of ways to improve the training of PINNs.
{The estimation error (c) 
comes from the use of finite data,
and is one of two components (together with approximation error) that constitute 
generalization error \cite{niyogi1999generalization}.
In machine learning, generalization error refers to a measure of the accuracy of the prediction on unseen data \cite{Mohri_18_FoundationsML}.
In PDE problems, the generalization error is
the distance between a global minimizer of the loss and the solution to the PDE,
where the distance has be to defined appropriately to reflect its regularity.
From this perspective, the present paper 
is concerned with the convergence of 
generalization error.}

{\bf Contributions.} In this work, we 
provide a convergence theory for PINNs 
with respect to the number of training data.
By adopting probabilistic space filling arguments \cite{Finlay_18_Lipschitz,Calder_19_consistency}, 
we derive an upper bound of the expected unregularized PINN loss \cite{Raissi_19_PINNs} (Theorem~\ref{thm:gen}).
Motivated by the upper bound,
we consider the problem of minimizing 
a specific regularized loss. 
By focusing on two classes of PDEs -- linear second order elliptic and parabolic that admit classical solutions -- 
we show that the sequence of minimizers of the regularized loss converges to the solution to the PDE uniformly. 
{In other words, generalization error of PINNs converges to zero under the uniform topology.}
In addition, we show that if minimizers satisfy the initial/boundary conditions, 
the mode of convergence becomes $H^1$.
To ensure the existence, regularity and uniqueness of the solution,
we adopt the Schauder approach.
To the best of our knowledge, this is the first theoretical work that proves the consistency of PINNs in the sample limit.
Computational examples are provided to demonstrate our theoretical findings.

The rest of the paper is organized as follows.
Upon briefly introducing mathematical setup in section~\ref{sec:setup}, 
we present the convergence analysis in section~\ref{sec:analysis}.
Computational examples are provided in section~\ref{sec:example}.
Finally, we conclude the paper in section~\ref{sec:conclusion}.

%% file: setup.tex
\section{Mathematical Setup and Preliminaries} \label{sec:setup}
Let $U$ be a bounded domain (open and connected) in $\mathbb{R}^d$.
We consider partial differential equations (PDEs) of the form 
\begin{equation} \label{def:governing}
\mathcal{L}[u](\x) =  f(\x) \quad \forall \x \in U, \qquad \mathcal{B}[u](\x) = g(\x) \quad \forall \x \in \Gamma \subseteq \partial U,
\end{equation}
where 
$\mathcal{L}[\cdot]$ is a differential operator and
$\mathcal{B}[\cdot]$ could be Dirichlet, Neumann, Robin, or periodic boundary conditions.
This paper considers PDEs that admit a unique classical solution $u(\x)$ with $\x=(x_1,\cdots,x_d)$.
The classical solution should satisfy the governing equation everywhere on $U$ and $\Gamma$.
We remark that for time-dependent problems, we consider time $t$ as a special component of $\x$, and $U$ contains the temporal domain. 
The initial condition can be simply treated as a special type of Dirichlet boundary condition on the spatio-temporal domain.

The goal is to approximate the solution to the PDE \eqref{def:governing}
from a set of training data. 
The training data consist of two types of data sets: residual and initial/boundary data.
A residual datum is a pair of input and output $(\x_r, f(\x_r))$,
where $\x_r \in U$
and an initial/boundary datum is a pair of input and output $(\x_{b}, g(\x_{b}))$, where $\x_{b} \in \Gamma$.
The set of $m_r$ residual input data points and the set of $m_{b}$ initial/boundary input data points are denoted by
$\mathcal{T}_r^{m_r} = \{\x_r^i\}_{i=1}^{m_r}$ 
and $\mathcal{T}_{b}^{m_{b}} = \{\x_{b}^i\}_{i=1}^{m_{b}}$, respectively.
Let us denote the vector of the number of training samples by $\bm{m}=(m_r,m_{b})$.
Note that we slightly abuse notation as in \cite{Raissi_19_PINNs}: $\x_r$ refers to a point in $U$.
Similarly, $\x_{b}$ refers to a point in $\Gamma$.
$m_r$ and $m_{b}$ represent the number of training data points in $U$ and $\Gamma$, respectively. 
\begin{remark}
    The present paper only considers the high regularity setting for PDEs 
    where point-wise evaluations 
    are well-defined.
    For the sake of simplicity, we consider one boundary operator case.
    And $u$ is a scalar real-valued function.
\end{remark}

Given a class of neural networks $\mathcal{H}_{n}$ 
that may depend on the number of training samples $\bm{m}$ (also may implicitly depend on the data sample itself),
we seek to find a neural network $h^*$ in $\mathcal{H}_{n}$ that minimizes 
an objective (loss) function.
Here $n$ represents a complexity of the class, e.g., the number of parameters.
To define an appropriate objective function, 
let us consider a loss integrand, which is the standard choice in physics informed neural networks (PINNs) \cite{Raissi_19_PINNs}:
\begin{equation} \label{def:loss}
\begin{split}
\textbf{L}(\x_r,{\x}_b;h,\bm{\lambda},\bm{\lambda}^R) &= \left(\lambda_r\|\mathcal{L}[h](\x_r)- f(\x_r)\|^2 \right)\mathbb{I}_{U}(\x_r) + \lambda_r^R R_r(h) \\
&+ 
\lambda_{b}\|\mathcal{B}[h](\x_{b}) - g(\x_{b})\|^2 \mathbb{I}_{\Gamma}(\x_{b})
+ \lambda_{b}^R R_{b}(h),
\end{split}
\end{equation}
where $\|\cdot\|$ is the Euclidean norm, 
$\mathbb{I}_{A}(\x)$ is the indicator function on the set $A$,
$\bm{\lambda} = (\lambda_r, \lambda_{b})$, 
$\bm{\lambda}^R = (\lambda_r^R, \lambda_{b}^R)$,
and $R_r(\cdot), R_{b}(\cdot)$ are regularization functionals.
Here $\lambda_r, \lambda_r^R,\lambda_{b}, \lambda_{b}^R \in \mathbb{R}_{+}\cup\{0\}$.

Suppose $\mathcal{T}_r^{m_r}$ and $\mathcal{T}_{b}^{m_{b}}$ are independently and identically distributed (iid) samples from probability distributions $\mu_r$ and $\mu_{b}$, respectively.
Let us define the empirical probability distribution on $\mathcal{T}_r^{m_r}$ 
by
$\mu_r^{m_r} = \frac{1}{m_r}\sum_{i=1}^{m_r} \delta_{\x_r^i}$.
Similarly, $\mu_{b}^{m_{b}}$ is defined.
The empirical loss and the expected loss are obtained by taking the expectations 
on the loss integrand \eqref{def:loss} with respect to $\mu^{\bm{m}}=\mu_r^{m_r}\times \mu_{b}^{m_{b}}$ and 
$\mu = \mu_r \times \mu_{b}$, respectively:
\begin{equation} \label{def:empirical-expected-loss}
	\begin{split}
	\text{Loss}_{\bm{m}}(h;\bm{\lambda},\bm{\lambda}^R) = \mathbb{E}_{\mu^{\bm{m}}}[\textbf{L}(\x_r,\x_b;h,\bm{\lambda},\bm{\lambda}^R)], 
	\quad
	\text{Loss}(h;\bm{\lambda},\bm{\lambda}^R) = \mathbb{E}_{\mu}[\textbf{L}(\x_r,\x_b;h,\bm{\lambda},\bm{\lambda}^R)].
	\end{split}
\end{equation}
In order for the expected loss to be well-defined,
it is assumed that $\mathcal{L}[h]$ and $f$ are in $L^2(U;\mu_r)$,
and $\mathcal{B}[h]$ and $g$ are in $L^2(\Gamma;\mu_{b})$
for all $h\in \mathcal{H}_{n}$.
If the expected loss function were available, 
its minimizer would be the solution to the PDE \eqref{def:governing} or close to it.
However, since it is unavailable in practice,
the empirical loss function is employed.
This leads to the following minimization problem:
\begin{equation} \label{def:problem}
\min_{h \in \mathcal{H}_{n}} \text{Loss}_{\bm{m}}(h;\bm{\lambda},\bm{\lambda}^R).
\end{equation}
We then hope a minimizer of the empirical loss to be close to the solution to the PDE \eqref{def:governing}.

We remark that in general, global minimizers to the problem of \eqref{def:problem} need not exist.
However, for $\epsilon > 0$, there always exists a 
$\epsilon$-suboptimal global solution $h^{\epsilon} \in \mathcal{H}_{n}$ \cite{Houska_19_Global}
satisfying $\text{Loss}_{\bm{m}}(h^{\epsilon}) \le \inf_{h \in \mathcal{H}_{n}} \text{Loss}_{\bm{m}}(h) + \epsilon$.
All the results of the present paper remain valid
if one replace global minimizers 
to $\epsilon$-suboptimal global minimizers for sufficiently small $\epsilon$.
Henceforth, for the sake of readability, we assume the existence of at least one global minimizer of the minimization problem \eqref{def:problem}.

When $\bm{\lambda}^R = 0$, we refer to the (either empirical or expected) loss as the (empirical or expected) PINN loss \cite{Raissi_19_PINNs}. 
We denote the empirical and the expected PINN losses by $\text{Loss}^{\text{PINN}}_{\bm{m}}(h;\bm{\lambda})$
and 
$\text{Loss}^{\text{PINN}}(h;\bm{\lambda})$, respectively.
Specifically,
\begin{equation} \label{def:PINN-loss}
\begin{split}
\text{Loss}^{\text{PINN}}_{\bm{m}}(h;\bm{\lambda})
&= \frac{\lambda_r}{m_r}\sum_{i=1}^{m_r} \|\mathcal{L}[h](\x_r^i)- f(\x_r^i)\|^2
+ \frac{\lambda_{b}}{m_{b}}\sum_{i=1}^{m_{b}} 
\|\mathcal{B}[h](\x_{b}^i) - g(\x_{b}^i)\|^2, \\
\text{Loss}^{\text{PINN}}(h;\bm{\lambda}) &= \lambda_r  \|\mathcal{L}[h]- f\|^2_{L^2(U;\mu)}
+\lambda_{b} \|\mathcal{B}[h] - g\|^2_{L^2(\Gamma;\mu_{b})}.
\end{split}
\end{equation}

Also, note that for $\bm{\lambda} \le \bm{\lambda'}$ 
and
$\bm{\lambda}^R \le \bm{\lambda'}^{R}$ 
(element-wise inequality), 
we have
$\text{Loss}_{\bm{m}}(h;\bm{\lambda},\bm{\lambda}^R) \le \text{Loss}_{\bm{m}}(h;\bm{\lambda}',\bm{\lambda'}^{R})$.


\subsection{Function Spaces and Regular Boundary}
Throughout this paper, we adopt the notation from \cite{Gilbarg_15_EllipticPDEs,Friedman_08_ParabolicPDEs}.
Let $U$ be a bounded domain in $\mathbb{R}^d$. 
Let $\x=(x_1,\cdots,x_d)$ be a point in $\mathbb{R}^d$.
For a positive integer $k$, 
let $C^{k}(U)$ be the set of functions having all derivatives of order $\le k$ continuous in $U$. 
Also, let $C^{k}(\overline{U})$ be the set of functions in $C^{k}(U)$ whose derivatives of order $\le k$ have continuous extensions to $\overline{U}$ (the closure of $U$).

We call a function $u$ uniformly H\"{o}lder continuous with exponent $\alpha$ in $U$ if the quantity
\begin{equation} \label{def:holder-coefficient}
[u]_{\alpha;U}=\sup_{x,y \in U, x\ne y} \frac{|u(x) - u(y)|}{\|x-y\|^\alpha} < \infty, \qquad 0 < \alpha \le 1,
\end{equation}
is finite.
Also, we call a function $u$ locally H\"{o}lder continuous with exponent $\alpha$ in $U$ if $u$ is uniformly H\"{o}lder continuous with exponent $\alpha$ on compact subsets of $U$. $[u]_{\alpha;U}$ is called the H\"{o}lder constant (coefficient) of $u$ on $U$.

Given a multi-index $\textbf{k}=(k_1,\cdots,k_d)$, we define
\begin{align*}
D^\textbf{k} u = \frac{\partial^{|\textbf{k}|}u}{\partial x_1^{k_1} \cdots \partial x_d^{k_d}},
\end{align*}
where $|\textbf{k}| = \sum_{i=1}^d k_i$.
Let $\textbf{k}=(k_1,\cdots,k_d)$ and $\textbf{k}'=(k'_1,\cdots,k'_d)$.
If $k_j' \le k_j$ for all $j$, we write $\textbf{k}' \le \textbf{k}$.
Let us define
\begin{align*}
[u]_{j,0;U} &:= \sup_{|\textbf{k}|=j} \sup_{U} |D^{\textbf{k}}u|, \qquad j=0,1, 2\cdots, \\
[u]_{j,\alpha;U} &:= \sup_{|\textbf{k}|=j} [D^{\textbf{k}}u]_{\alpha;U} =  \sup_{|\textbf{k}|=j} \left[\sup_{x,y \in U, x\ne y}  \frac{\|D^{\textbf{k}}u(x) - D^{\textbf{k}}u(y)\|}{\|x-y\|^\alpha}\right].
\end{align*}

\begin{definition}
	For a positive integer $k$
	and $0 < \alpha \le 1$, 
	the H\"{o}lder spaces $C^{k,\alpha}(\overline{U})$ ($C^{k,\alpha}(U)$) are 
	the subspaces of $C^{k}(\overline{U})$ ($C^{k}(U)$) consisting of 
	all functions $u \in C^{k}(\overline{U})$ ($C^{k}(U)$)
	satisfying $\sum_{j=0}^{k} [u]_{j,0;U} + [u]_{k,\alpha;U} < \infty$.
\end{definition}
For simplicity, we often write $C^{k,0} = C^{k}$ and $C^{0,\alpha} = C^\alpha$ for $0 < \alpha < 1$.
The related norms are defined on $C^{k}(\overline{U})$ and $C^{k,\alpha}(\overline{U})$, respectively, by 
\begin{align*}
\|u\|_{C^{k}(\overline{U})} =  \sum_{j=0}^{k} [u]_{j,0;U}, \qquad
\|u\|_{C^{k,\alpha}(\overline{U})} = \|u\|_{C^{k}(\overline{U})} + [u]_{k,\alpha;U}, \qquad 0 < \alpha \le 1.
\end{align*}
With these norms, $C^{k}(\overline{U})$ and $C^{k,\alpha}(\overline{U})$ are Banach spaces. 
Also, we denote $\frac{\partial u}{\partial x_j}$ as $D_ju$
and $\frac{\partial^2 u}{\partial x_i \partial x_j}$ as $D_{ij}u$.

\begin{definition}
	A bounded domain $U$ in $\mathbb{R}^d$ and its boundary are said to be of class $C^{k,\alpha}$ where $0\le \alpha \le 1$, if
	at each point $x_0 \in \partial U$
	there is a ball $B=B(x_0)$ and a one-to-one mapping
	$\psi$ of $B$ onto $D \subset \mathbb{R}^d$ such that 
	(i) $\psi(B\cap U) \subset \mathbb{R}_{+}^d$,
	(ii) $\psi(B\cap \partial U) \subset \partial \mathbb{R}_{+}^d$,
	(iii) 
	$\psi \in C^{k,\alpha}(B), \psi^{-1} \in C^{k,\alpha}(D)$.
\end{definition}

For the parabolic equations, 
we introduce function spaces for interior estimates
following \cite{Friedman_08_ParabolicPDEs}.
For $T > 0$ and a bounded domain $U$ in $\mathbb{R}^d$, 
let $\Omega = U \times (0,T)$,
$B = U \times \{0\}$, $S_{\tau} = \partial U \times [0, \tau]$.
For any two points $P=(\x,t)$ and $Q=(\x', t')$ in $\Omega$,
let $d(P,Q) = (\|\x-\x'\|^2 + |t-t'|)^{1/2}$
where $\|\cdot\|$ is the Euclidean norm.
For $Q = (\x, \tau)$, 
let $d_Q:= \inf_{P \in B\cup S_{\tau}} d(P,Q)$ be the distance from $Q$ to $B\cup S_{\tau}$.
For any two points in $P, Q \in \Omega$,
let $d_{PQ} = \min(d_P, d_Q)$.
Let $D_t =\partial/\partial t$.
Let us define 
\begin{equation} \label{def:parabolic-norm}
    \begin{split}
        |u|_{\alpha} &= \sup_{P \in \Omega} |u(P)| + \sup_{P, Q \in \Omega} d_{PQ}^\alpha \frac{|u(P) - u(Q)|}{d(P,Q)^\alpha}, \\
        |\text{d}^m u|_{\alpha} &= \sup_{P \in \Omega} d_P^m |u(P)| + \sup_{P,Q \in \Omega} d_{PQ}^{m+\alpha} \frac{|u(P) - u(Q)|}{d(P,Q)^\alpha}, \\
        |u|_{2+\alpha} &= |u|_{\alpha} + \sum_{i} |\text{d}D_i u|_{\alpha}
    + \sum_{i,j}|\text{d}D_{ij} u|_{\alpha} + |\text{d}^2D_t u|_{\alpha}.
    \end{split}
\end{equation}
Let $C_{\alpha}(\Omega)$ and $C_{2+\alpha}(\Omega)$ be the space of all functions $u$
with finite norm $|u|_{\alpha}$ and $|u|_{2+\alpha}$,
respectively. 
They are Banach spaces.
We refer the readers to \cite{Friedman_08_ParabolicPDEs}
for more details.

\subsection{Neural Networks}
Let $h^L: \mathbb{R}^{d} \to \mathbb{R}^{d_{\text{out}}}$ 
be a feed-forward neural network having $L$ layers and $n_\ell$ neurons in the $\ell$-th layer. 
The weights and biases in the $l$-th layer are represented by 
a weight matrix $\mathbf{W}^l \in \mathbb{R}^{n_l \times n_{l-1}}$ 
and a bias vector $\mathbf{b}^l \in \mathbb{R}^{n_l}$, respectively. 
Let $\bm{\theta}_L:=\{\bm{W}^j, \bm{b}^j\}_{1\le j \le L}$.
For notational completeness, let $n_0 = d$ and $n_L = d_{\text{out}}$.
For a fixed positive integer $L$, let $\vec{\bm{n}} = (n_0,n_1,\cdots,n_L) \in \mathbb{N}^{L+1}$ where $\mathbb{N}=\{1,2,3, \cdots \}$.
Then, $\vec{\bm{n}}$ describes a network architecture.
Given an activation function $\sigma(\cdot)$ which is applied element-wise,
the feed-forward neural network is defined by
\begin{align*}
	h^{\ell}(\textbf{x}) &= \bm{W}^{\ell}\sigma(h^{\ell-1}(\textbf{x})) + \bm{b}^{\ell}
	\in \mathbb{R}^{N_\ell}, 
	\qquad \text{for} \quad 2 \le \ell \le L
\end{align*}
and $h^1(\textbf{x}) = \bm{W}^{1}\textbf{x} + \bm{b}^{1}$.
The input is $\mathbf{x} \in \mathbb{R}^{n_0}$, 
and the output of the $\ell$-th layer is $h^\ell(\x) \in \mathbb{R}^{n_\ell}$. 
Popular choices of activation functions include 
the sigmoid ($1/(1+e^{-x})$), the hyperbolic tangent ($\tanh(x)$),
and the rectified linear unit ($\max\{x,0\}$).
Note that $h^L$ is called a $(L-1)$-hidden layer neural network or a $L$-layer neural network.

Since a network $h^L(\x)$ depends on the network parameters $\bm{\theta}_L$
and the architecture $\vec{\bm{n}}$,
we often denote $h^L(\x)$ by $h^L(\x;\vec{\bm{n}}, \bm{\theta}_L)$. 
If $\vec{\bm{n}}$ is clear in the context, we simple write $h^L(\x;\bm{\theta}_L)$.
Given a network architecture, we define a neural network function class
\begin{equation} \label{def:NN-class}
	\mathcal{H}_{\vec{\bm{n}}}^{\text{NN}} = \left\{h^L(\cdot;\vec{\bm{n}},\bm{\theta}_L):\mathbb{R}^{d}\mapsto \mathbb{R}^{d_\text{out}} | \bm{\theta}_L = \{(\bm{W}^j,\bm{b}^j) \}_{j=1}^L
	  \right\}.
\end{equation}

Since $h^L(\x;\bm{\theta}_L)$ is parameterized by $\bm{\theta}_L$, 
the problem of \eqref{def:problem} with $\mathcal{H}_{\vec{\bm{n}}}^{\text{NN}}$ is
equivalent to
\begin{equation*} 
	\min_{\bm{\theta}_L}  \text{Loss}_{\bm{m}}(\bm{\theta}_L;\bm{\lambda},\bm{\lambda}^R), \quad
	\text{where} \quad
	\text{Loss}_{\bm{m}}(\bm{\theta}_L;\bm{\lambda},\bm{\lambda}^R) = \text{Loss}_{\bm{m}}(h(\x;\bm{\theta}_L);\bm{\lambda},\bm{\lambda}^R).
\end{equation*}
Throughout this paper, the activation function is assumed to be sufficiently smooth,
which is common in practice.

\begin{lemma} \label{lem:NN}
	Let $U$ be a bounded domain
	and $\mathcal{H}_{\vec{\bm{n}}}^{NN}$ be a class of neural networks whose architecture is 
	$\vec{\bm{n}} = (n_0,\cdots,n_L)$
	whose activation function $\sigma(x) \in C^{k'}(\mathbb{R})$
	satisfies that for each $s\in \{0,\cdots,k\}$ where $k < k'$,  
	$\frac{d^s \sigma(x)}{dx^s}$ is bounded and Lipschitz continuous.
	For a function $u \in C^k(\overline{U})$,
	let $\{h_j\}_{j=1}^\infty$ be a sequence of networks in $\mathcal{H}_{\vec{\bm{n}}}^{NN}$ such that 
	the associated weights and biases are uniformly bounded
	and 
	$h_j \to u$ in $C^0(\overline{U})$.
	Then, $h_j \to u$ in $C^k(\overline{U})$.
\end{lemma}
\begin{proof}
	The proof can be found in Appendix~\ref{app:lem:NN}.
\end{proof}

It can be checked that the $\tanh$ activation function satisfies the conditions of Lemma~\ref{lem:NN} for all $k$ as follows:
Note that $\tanh(x)$ is bounded by 1 and is 1-Lipschitz continuous. 
The $k$-th derivative of $\tanh(x)$
is expressed as a polynomial of $\tanh(x)$ with a finite degree,
which shows both the boundedness and Lipschitz continuity.

\begin{remark}
    In what follows, 
    we simply write 
    $\mathcal{H}_{\vec{\bm{n}}}^{NN}$
    as
    $\mathcal{H}_{\bm{m}}$, 
    assuming $\vec{\bm{n}}$ depends on $\bm{m}$
    and possibly implicitly on the data samples itself.
    Since neural networks
    are universal approximators,
    the network architecture $\vec{\bm{n}}$ is expected to 
    grow 
    proportionally on $\bm{m}$.
\end{remark}

%% file: analysis.tex
\section{Convergence Analysis} \label{sec:analysis}
If the expected loss function were available,
a function that minimizes it would be sought.
However, since the expected loss is not available in practice,
the empirical loss function is employed.
We first derive an upper bound of the expected PINN loss \eqref{def:PINN-loss}.
The bound involves a specific regularized empirical loss.

The derivation is based on the probabilistic space filling arguments \cite{Calder_19_consistency,Finlay_18_Lipschitz}.
In this regard, we make the following assumptions on the training data distributions.
\begin{assumption} \label{assumption:data-dist}
	Let $U$ be a bounded domain in $\mathbb{R}^d$
	that is at least of class $C^{0,1}$ 
	and $\Gamma$ be a closed subset of $\partial U$.
	Let $\mu_r$ and $\mu_{b}$ be probability distributions defined on $U$ and $\Gamma$, respectively.
	Let $\rho_r$ be the probability density of $\mu_r$ with respect to $d$-dimensional Lebesgue measure on $U$.
	Let $\rho_{b}$ be the probability density of $\mu_{b}$
	with respect to $(d-1)$-dimensional Hausdorff measure on $\Gamma$.
	\begin{enumerate}
		\item $\rho_r$ and $\rho_{b}$ are supported on $\overline{U}$
		and $\Gamma$, respectively.
		Also, $\inf_{U} \rho_r > 0$ and $\inf_{\Gamma} \rho_{b} > 0$.
		\item 
		For $\epsilon > 0$, there exists partitions of $U$ and $\Gamma$, $\{U_j^\epsilon\}_{j=1}^{K_{r}}$
		and $\{\Gamma_j^\epsilon\}_{j=1}^{K_{b}}$
		that depend on $\epsilon$
		such that 
		for each $j$, 
		there are cubes $H_{\epsilon}(\textbf{z}_{j,r})$ and 
		$H_{\epsilon}(\textbf{z}_{j,b})$ of side length $\epsilon$
		centered at $\textbf{z}_{j,r} \in U_j^\epsilon$
		and $\textbf{z}_{j,b} \in \Gamma_j^\epsilon$, respectively, 
		satisfying $U_j^\epsilon \subset H_{\epsilon}(\textbf{z}_{j,r})$
		and $\Gamma_j^\epsilon \subset H_{\epsilon}(\textbf{z}_{j,b})$.
		\item 
		There exists positive constants $c_r, c_{b}$ such that 
		$\forall \epsilon > 0$,
		the partitions from the above satisfy
		$c_r \epsilon^{d} \le \mu_r(U_j^\epsilon)$
		and
		$c_{b} \epsilon^{d-1} \le 
		\mu_{b}(\Gamma_j^\epsilon)$
		for all $j$.

		There exists positive constants $C_r, C_{b}$ such that 
		$\forall \x_r \in U$ and $\forall \x_b \in \Gamma$,
		$\mu_r(B_{\epsilon}(\x_r) \cap U) \le C_r\epsilon^d$
		and 
		$\mu_{b}(B_\epsilon(\x_{b}) \cap \Gamma) \le C_{b} \epsilon^{d-1}$
		where 
		$B_\epsilon(\x)$ is a closed ball of radius $\epsilon$ centered at $\x$.

		Here $C_r, c_r$ depend only on $(U,\mu_r)$
		and $C_{b}, c_{b}$ depend only on $(\Gamma, \mu_{b})$.
		\item When $d=1$, 
		we assume that all boundary points are available.
		Thus, no random sample is needed on the boundary.
	\end{enumerate}
\end{assumption}

We remark that Assumption~\ref{assumption:data-dist} guarantees that random samples drawn from probability distributions can fill up both the interior of the domain $U$ and the boundary $\partial U$. 
These are mild assumptions and can be satisfied in many practical cases. 
For example, let $U=(0,1)^d$. Then the uniform probability distributions on both $U$ and $\partial U$ satisfy Assumption~\ref{assumption:data-dist}.


We now state our result that bounds the expected PINN loss in terms of a regularized empirical loss.  
Let us recall that $\bm{m}$ is the vector of the number of training data points, i.e.,
$\bm{m} = (m_r,m_{b})$.
The constants $c_r, C_r, c_{b}, C_{b}$ are introduced in Assumption~\ref{assumption:data-dist}.
For a function $u$, $[u]_{\alpha;U}$ is the H\"{o}lder constant of $u$ with exponent $\alpha$ in $U$ \eqref{def:holder-coefficient}.
\begin{theorem} \label{thm:gen}
	Suppose Assumption~\ref{assumption:data-dist} holds.
	Let $m_r$ and $m_{b}$ be the number of iid samples from $\mu_r$ and $\mu_{b}$, respectively.
	For some $0 < \alpha \le 1$, 
	let $h, f, g$, $R_r(h)$, $R_{b}(h)$  satisfy
	\begin{align*}
	    \big[\mathcal{L}[h]\big]^2_{\alpha;U} \le R_r(h)  < \infty, \quad
	    \big[\mathcal{B}[h]\big]_{\alpha;\Gamma}^2 \le R_{b}(h) < \infty, \quad
	    \big[f\big]_{\alpha;U}, \big[g\big]_{\alpha;\Gamma}  < \infty.
	\end{align*}
	Let $\bm{\lambda} = (\lambda_r,\lambda_{b})$ be a fixed vector.
    Let $\bm{\hat{\lambda}}_{\bm{m}}^R = (\hat{\lambda}_{r,\bm{m}}^R,\hat{\lambda}_{b,\bm{m}}^R)$
    be a vector whose elements to be defined.

	For $d \ge 2$, with probability at least, 
	$(1 - \sqrt{m_r}(1-1/\sqrt{m_r})^{m_r})(1 - \sqrt{m_{b}}(1-1/\sqrt{m_{b}})^{m_{b}})$,
	we have
	\begin{equation*}
	    \text{Loss}^{\text{PINN}}(h;\bm{\lambda})
        \le C_{\bm{m}}\cdot \text{Loss}_{\bm{m}}(h;\bm{\lambda}, \bm{\hat{\lambda}}_{\bm{m}}^R)
        +C'(m_r^{-\frac{\alpha}{d}} + m_{b}^{-\frac{\alpha}{d-1}}),
	\end{equation*}
	where $\kappa_r = \frac{C_r}{c_r}$, $\kappa_{b} = \frac{C_{b}}{c_{b}}$, 
	$C_{\bm{m}} = 3\max\{\kappa_r \sqrt{d}^{d} m_r^{\frac{1}{2}}, \kappa_{b} \sqrt{d}^{d-1} m_{b}^{\frac{1}{2}}\}$, 
	$C'$ is a universal constant that depends only on $\bm{\lambda}$, $d$, $c_r$, $c_{b}$, $\alpha$, $f$, $g$,
	and
	\begin{equation} \label{def:lambdas}
	\begin{split}
	    \hat{\lambda}_{r,\bm{m}}^R = \frac{3\lambda_r \sqrt{d}^{2\alpha}c_r^{-\frac{2\alpha}{d}}}{C_{\bm{m}}}\cdot m_r^{-\frac{\alpha}{d}},
        \quad
        \hat{\lambda}_{b,\bm{m}}^R = \frac{3\lambda_{b} \sqrt{d}^{2\alpha} c_{b}^{-\frac{2\alpha}{d-1}}}{C_{\bm{m}}}\cdot m_{b}^{-\frac{\alpha}{d-1}}.
	\end{split}
	\end{equation}

	For $d = 1$,
	with probability at least, $1 - \sqrt{m_r}(1-1/\sqrt{m_r})^{m_r}$, 
	we have
	\begin{equation*}
	\text{Loss}^\text{PINN}(h;\bm{\lambda})
	\le C_{\bm{m}}\cdot \text{Loss}_{\bm{m}}(h;\bm{\lambda},\bm{\hat{\lambda}}^R_{\bm{m}}) +C'm_r^{-\alpha},
	\end{equation*}
	where 
	$C_{\bm{m}} = 3\kappa_r m_r^{\frac{1}{2}}$, 
	$\hat{\lambda}^R_{r,\bm{m}} = \frac{\lambda_r c_r^{-2\alpha}}{\kappa_r}\cdot m_r^{-\alpha-\frac{1}{2}}$,
	$\hat{\lambda}_{b,\bm{m}}^R = 0$ 
	and $C'$ is a universal constant that depends only on $\lambda_r$, $c_r$, $\alpha$, $f$.
\end{theorem}
\begin{proof}
	The proof can be found in Appendix~\ref{app:thm:gen}.
\end{proof}

Let $\bm{\lambda}$ be a vector independent of $\bm{m}$ 
and 
$\bm{\lambda}_{\bm{m}}^R = (\lambda_{r,\bm{m}}^R, \lambda_{b,\bm{m}}^R)$
be a vector satisfying
\begin{equation} \label{def:lambda-condition}
    \bm{\lambda}_{\bm{m}}^R \ge \bm{\hat{\lambda}}_{\bm{m}}^R, \qquad
    \|\bm{\lambda}_{\bm{m}}^R\|_\infty = \mathcal{O}(\|\bm{\hat{\lambda}}_{\bm{m}}^R\|_\infty),
\end{equation}
where $\bm{\hat{\lambda}}_{\bm{m}}^R$ is defined in \eqref{def:lambdas}.
For a vector $v$, $\|v\|_{\infty}$ is the maximum norm, i.e., $\|v\|_{\infty} = \max_i |v_i|$.
We note that 
since $\bm{\hat{\lambda}}_{\bm{m}}^R \to 0$
as $\bm{m} \to \infty$,
the above condition implies 
$\bm{\lambda}_{\bm{m}}^R \to 0$
as $\bm{m} \to \infty$.	

By letting $R_r(h) = \big[\mathcal{L}[h]\big]_{\alpha;U}^2$
and $R_{b}(h) = \big[\mathcal{B}[h]\big]_{\alpha;\Gamma}^2$,
let us define the H\"{o}lder regularized empirical loss:
\begin{equation} \label{def:Holder-Reg-Loss}
\begin{split}
&\text{Loss}_{\bm{m}}(h;\bm{{\lambda}},\bm{{\lambda}}^R_{\bm{m}}) 
\\
&= \begin{cases}
\text{Loss}_{\bm{m}}^{\text{PINN}}(h;\bm{{\lambda}})
+\lambda_{r,\bm{m}}^R
\big[\mathcal{L}[h]\big]_{\alpha;U}^2
+\lambda_{b,\bm{m}}^R
\big[\mathcal{B}[h]\big]_{\alpha;\Gamma}^2, & \text{if } d \ge 2 \\[12pt]
\text{Loss}_{\bm{m}}^{\text{PINN}}(h;\bm{{\lambda}})
+ \lambda_{r,\bm{m}}^R
\big[\mathcal{L}[h]\big]_{\alpha;U}^2,  &\text{if } d = 1
\end{cases}
\end{split}
\end{equation}
where $\bm{\lambda}_{\bm{m}}^R$ are vectors satisfying \eqref{def:lambda-condition}.
We note that the H\"{o}lder regularized loss \eqref{def:Holder-Reg-Loss} is greater than or equal to the empirical loss shown in Theorem~\ref{thm:gen} (assuming $R_r(h) = \big[\mathcal{L}[h]\big]_{\alpha;U}^2$
and $R_{b}(h) = \big[\mathcal{B}[h]\big]_{\alpha;\Gamma}^2$).
Since $\big[\mathcal{L}[h]\big]_{\alpha;U}^2$
and $\big[\mathcal{B}[h]\big]_{\alpha;\Gamma}^2$ do not depend on the training data
and $\hat{\lambda}_{r,\bm{m}}^R, \hat{\lambda}_{b,\bm{m}}^R \to 0$ as $m_r, m_{b} \to \infty$,
this suggests that the more data we have, the less regularization is needed.

Theorem~\ref{thm:gen} indicates that 
minimizing the H\"{o}lder regularized loss \eqref{def:Holder-Reg-Loss} over $h$
results in minimizing an upper bound of the expected PINN loss \eqref{def:PINN-loss}.
Thus the H\"{o}lder regularized loss could be used as a loss function in a way to have a small expected PINN loss.
In general, however, the H\"{o}lder-regularization functionals (the H\"{o}lder constants) are impractical to be evaluated numerically. Also, minimizing an upper bound does not necessarily imply minimizing the expected PINN loss.

If the domain is convex and the H\"{o}lder exponent is $\alpha = 1$ (i.e., Lipschitz constant),
it follows from Rademacher's Theorem \cite{Evans_15_measure} that 
the Lipschitz constant of $\mathcal{L}[h]$
is the supremum of the sup norm of $\nabla \mathcal{L}[h]$ over $U$ (assuming $\nabla \mathcal{L}[h]$ exists).
Thus, in practice, one can use the maximum of the sup norm of the derivative over the set of training data points
to estimate the Lipschitz constant.
The resulting Lipschitz regularized (LIPR) loss is given by
\begin{equation} \label{def:LIPR-loss}
\begin{split}
\text{Loss}_{\bm{m}}^{\text{LIPR}}(h;\bm{\lambda},\bm{\lambda}^R) &= 
\text{Loss}_{\bm{m}}^{\text{PINN}}(h;\bm{\lambda}) +
\lambda_{r}^R
\max_{1 \le i \le m_r} \|\nabla \mathcal{L}[h](\x_r^i)\|_{\infty}^2
+
\lambda_{b}^R 
\max_{1\le j\le m_{b}} \|\nabla \mathcal{B}[h](\x_{b}^j)\|_{\infty}^2,
\end{split}
\end{equation}
where $\text{Loss}_{\bm{m}}^{\text{PINN}}(h;\bm{\lambda})$ is defined in \eqref{def:PINN-loss}
and $\bm{\lambda}^R = (\lambda_r^R, \lambda_{b}^R)$ are weights of the regularization penalty terms.

\subsection{Expected PINN Loss} 
We now consider the problem of minimizing the H\"{o}lder regularized loss \eqref{def:Holder-Reg-Loss}.
In particular, we want to quantify how well a minimizer of the regularized empirical loss
performs on the expected PINN loss \eqref{def:PINN-loss}.

We make the following assumptions 
on the classes of neural networks
for the minimization problems \eqref{def:problem}.
\begin{assumption} \label{assumption:convergence}
	Let $k$ be the highest order of the derivative shown in the PDE \eqref{def:governing}.
	For some $0 < \alpha \le 1$, 
	let $f \in C^{0,\alpha}(U)$ and $g \in C^{0,\alpha}(\Gamma)$.
	\begin{enumerate}
	    \item For each $\bm{m}$,
	    let $\mathcal{H}_{\bm{m}}$ be a class of neural networks
	    in $C^{k,\alpha}(U)\cap C^{0,\alpha}(\overline{U})$
	    such that
    	for any $h \in \mathcal{H}_{\bm{m}}$, $\mathcal{L}[h] \in C^{0,\alpha}(U)$ 
	    and $\mathcal{B}[h] \in C^{0,\alpha}(\Gamma)$.
		\item For each $\bm{m}$, $\mathcal{H}_{\bm{m}}$ contains a network $u_{\bm{m}}^*$ 
		satisfying 
		$\text{Loss}_{\bm{m}}^{\text{PINN}}(u_{\bm{m}}^*;\bm{\lambda}) = 0$.
		\item 
		And, 
		$$
		\sup_{\bm{m}} \big[\mathcal{L}[u_{\bm{m}}^*]\big]_{\alpha;U}<\infty, \quad \sup_{\bm{m}} \big[\mathcal{B}[u_{\bm{m}}^*]\big]_{\alpha;\Gamma} < \infty.
		$$
	\end{enumerate}
\end{assumption}

All the assumptions are essential for the proof.
Our proof relies on 
the use of 
Assumption~\ref{assumption:convergence} 
that guarantees the uniform equicontinuity of subsequences.
All assumptions hold automatically if $\mathcal{H}_{\bm{m}}$ contains the solution to the PDE for all $\bm{m}$.
For example, \cite{Darbon_19_NNsHJ,Darbon_20_NNsHJ} showed that the solution 
to some Hamilton-Jacobi PDEs can be exactly represented by neural networks.
The second assumption could be relaxed, however, we do not discuss it here for the readability (See Appendix~\ref{app:thm:conv-loss}).
{
\begin{remark}
    The choice of classes of neural networks $\mathcal{H}_{\bm{m}}$
    may depend on many factors
    including 
    the underlying PDEs, 
    the training data,
    and the network architecture.
    Assumption~\ref{assumption:convergence}
    provides a set of conditions 
    for neural networks 
    for the purpose of our analysis.
    We will not discuss this issue further in the present paper 
    as it beyond our scope. 
\end{remark}}

For the rest of this paper, we use the following notation.
When the number of the initial/boundary training data points $m_{b}$ is completely determined by the number of residual points $m_r$ (e.g. $m_r^{d-1} = \mathcal{O}(m_{b}^{d})$),
the vector of the number of training data $\bm{m}$ depends only on $m_r$. 
In this case, we simply write $\mathcal{H}_{\bm{m}}$, $\bm{\lambda}_{\bm{m}}^R$, $\text{Loss}_{\bm{m}}$ 
as $\mathcal{H}_{m_r}$, $\bm{\lambda}_{m_r}^R$, $\text{Loss}_{m_r}$, respectively.

We now show that minimizers of the H\"{o}lder regularized loss \eqref{def:Holder-Reg-Loss} indeed produce a small expected PINN loss.
\begin{theorem} \label{thm:conv-loss}
	Suppose Assumptions~\ref{assumption:data-dist} and ~\ref{assumption:convergence} hold.
	Let $m_r$ and $m_{b}$ be the number of iid samples from 
	$\mu_r$ and $\mu_{b}$, respectively,
	and $m_r = \mathcal{O}(m_{b}^{\frac{d}{d-1}})$.
	Let $\bm{\lambda}_{m_r}^R$ be a vector satisfying 
	\eqref{def:lambda-condition}.
	Let $h_{m_r} \in \mathcal{H}_{m_r}$ be a minimizer of the H\"{o}lder regularized loss $\text{Loss}_{m_r}(\cdot;\bm{\lambda}, \bm{\lambda}_{m_r}^R)$ \eqref{def:Holder-Reg-Loss}.
	Then the following holds.
	\begin{itemize}
	    \item With probability at least $(1 - \sqrt{m_r}(1-c_r/\sqrt{m_r})^{m_r})(1 - \sqrt{m_{b}}(1-c_{b}/\sqrt{m_{b}})^{m_{b}})$
	    over iid samples,
	\begin{equation*}
	\text{Loss}^{\text{PINN}}(h_{m_r};\bm{\lambda}) = \mathcal{O}(m_r^{-\frac{\alpha}{d}}).
	\end{equation*} 
	    \item With probability 1 over iid samples, 
	\begin{equation} \label{thm:l2-convergence}
	\lim_{m_r \to \infty} \mathcal{L}[h_{m_r}] = f \text{ in } C^0(U),
	\quad
	\lim_{m_r \to \infty} \mathcal{B}[h_{m_r}] = g \text{ in } C^0(\Gamma).
	\end{equation}
	\end{itemize}
\end{theorem}
\begin{proof}
	The proof can be found in Appendix~\ref{app:thm:conv-loss}.
\end{proof}

\begin{remark}
    {The main goal of this section is to 
establish 
the convergence of $\mathcal{L}[h_{m_r}]$ and $\mathcal{B}[h_{m_r}]$ as ${m_r} \to \infty$.
Theorem~\ref{thm:conv-loss} provides 
a general convergence result under the
aforementioned conditions
on classes of neural networks.
The derived rate of convergence
is largely due to the probabilistic space filling argument
\cite{Finlay_18_Lipschitz,Calder_19_consistency}.
}
\end{remark}

Theorem~\ref{thm:conv-loss} shows that the expected PINN loss \eqref{def:PINN-loss} at minimizers of the H\"{o}lder regularized losses \eqref{def:Holder-Reg-Loss} converges to zero.
And more importantly, it shows the uniform convergences of $\mathcal{L}[h_{m_r}] \to f$ and $\mathcal{B}[h_{m_r}] \to g$ as ${m_r} \to \infty$.
(In fact, the mode of convergence is $C^{0,\beta}$ for all
$\beta < \alpha$ by the compact embedding of H\"{o}lder spaces).
These results are meaningful, however,
it is insufficient to claim the convergence of $h_{m_r}$ to the solution to the PDE \eqref{def:governing}.
In the next subsections, we consider two classes of PDEs and discuss the convergence of neural networks.

Remark that the uniform convergences of \eqref{thm:l2-convergence}
do not necessarily imply the convergence of $h_m$ to the solution $u^*$. 
For example, let $f_n(x) = n\cos(x/n)$ on $U=(-1,1)$.
It can be checked that for any $k \ge 1$, the $k$-th derivative of $f_n$ converges to $0$ uniformly.
However, $f_n(x)$ does not even converge.

\subsection{Linear Elliptic PDEs}
We consider the second-order linear elliptic PDEs
with the Dirichlet boundary condition:
\begin{equation} \label{def:elliptic-PDE}
\begin{split}
\begin{cases}
\mathcal{L}[u] = f  & \text{ in } U, \\
u = g      & \text{ in } \partial U,
\end{cases}
\end{split}
\end{equation}
where $f:U \to \mathbb{R}$, $g:\partial U \mapsto \mathbb{R}$, and 
\begin{align*}
\mathcal{L}[u] = \sum_{i,j=1}^{d} D_i\left(a^{ij}(\x) D_ju + b^i(x)u\right) + \sum_{i=1}^{d} c^i(\x) D_iu + d(\x)u.
\end{align*} 
The coefficients are defined on $U$.
To guarantee the existence, regularity and uniqueness of the classical solution to  \eqref{def:elliptic-PDE}, 
we adopt the Schauder approach (Chapter 6 of \cite{Gilbarg_15_EllipticPDEs}).
For the proof, the following assumptions are made.
 
\begin{assumption} \label{assumption:elliptic}
	Let $\lambda(\x)$ be the minimum eigenvalues of $[a^{ij}(\x)]$
	and $\alpha \in (0,1)$.
	\begin{enumerate}
		\item (Uniformly elliptic) For some constant $\lambda_0$, $\lambda(\x) \ge \lambda_0 > 0$ in $U$
		and $a^{ij} = a^{ji}$.
		\item The coefficients $c^i, d$ of the operator $\mathcal{L}$ are in $C^{\alpha}(U)$.
		\item $a^{ij}$ and $b^i$ are in $C^{1, \alpha}(U)$.
		Also, 
		$f \in C^{0,\alpha}(U)$
	    and $g \in C^{0,\alpha}(\partial U)$, 
		\item $U$ satisfies the exterior sphere condition at every boundary point.
		\item There are constants $\Lambda, v \ge 0$ such that for all $x \in U$
		\begin{equation*}
		    \sum_{i,j} |a^{ij}(\x)|^2 \le \Lambda, \qquad \lambda_0^{-2}(\|\bm{b}(\x)\|^2 + \|\bm{c}(\x)\|^2) + \lambda_0^{-1}|d(\x)| \le v^2,
		\end{equation*}
		where $\bm{b}(\x) = [b^1(\x),\cdots,b^d(\x)]^T$
		and $\bm{c}(\x) = [c^1(\x),\cdots, c^d(\x)]^T$.
	\end{enumerate}
\end{assumption}



We now present our main convergence result for the linear elliptic PDEs.
\begin{theorem} \label{thm:main-elliptic}
	Suppose Assumptions~\ref{assumption:data-dist}, ~\ref{assumption:convergence}
	and ~\ref{assumption:elliptic} hold
	and $\tilde{c}(\x)  = \sum_{i=1}^{d} D_ib^i(\x)+d(\x) \le 0$.
	Let $m_r$ and $m_{b}$ be the number of iid samples from 
	$\mu_r$ and $\mu_{b}$, respectively,
	and $m_r = \mathcal{O}(m_{b}^{\frac{d}{d-1}})$.
	Let $\bm{\lambda}_{m_r}^R$ be a vector satisfying 
	\eqref{def:lambda-condition}.
	Let $h_{m_r} \in \mathcal{H}_{m_r}$ be a minimizer of the H\"{o}lder regularized loss $\text{Loss}_{m_r}(\cdot;\bm{\lambda}, \bm{\lambda}_{m_r}^R)$ \eqref{def:Holder-Reg-Loss}.
	Then the following holds.
	\begin{enumerate}
		\item For each PDE data set $\{f, g\}$ satisfying 
		Assumption~\ref{assumption:elliptic},
	    there exists a unique classical solution $u^*$ to the PDE \eqref{def:elliptic-PDE}.
		\item With probability 1 over iid samples, 
		\begin{equation*}
		    \lim_{m_r \to \infty} h_{m_r} = u^*, \qquad \text{in } C^0(U).
		\end{equation*}
	\end{enumerate}
\end{theorem}
\begin{proof}
	The proof can be found in Appendix~\ref{app:thm:main-elliptic}.
\end{proof}

{
Theorem~\ref{thm:main-elliptic} shows that neural networks that minimize the H\"{o}lder regularized empirical losses \eqref{def:Holder-Reg-Loss}
converge to the unique classical solution to the PDE \eqref{def:elliptic-PDE}.
This answers the question
we posed in the introduction
for linear second-order elliptic PDEs.
It shows the consistency of
PINNs in the sample limit.
Equivalently, 
it can be interpreted as 
the convergence of generalization error
measured with the uniform topology.
}

Next, we show that if each minimizer exactly satisfies the boundary conditions, 
the mode of convergence becomes $H^1$.

\begin{theorem} \label{thm:elliptic-H1}
	Under the same conditions of Theorem~\ref{thm:main-elliptic}, 
	suppose each minimizer of \eqref{def:problem}
	satisfies the boundary conditions.
	Then, with probability 1 over iid samples, 
	the sequence of minimizers stated in Theorem~\ref{thm:main-elliptic}
	converges to the solution $u^*$ to the PDE \eqref{def:elliptic-PDE} in $H^1(U)$.
\end{theorem}
\begin{proof}
	The proof can be found in Appendix~\ref{app:thm:elliptic-H1}.
\end{proof}

In the literature, there are several ways to enforce neural networks to satisfy the boundary conditions. 
The work of \cite{Lagaris_98_ANN-ODE-PDE} considered the function classes that exactly satisfy the boundary conditions. The idea was then extended and generalized in \cite{Lagaris_00_ANN-Irregular,Berg_18_Unified}
to deal with irregular domains. 
This approach consists of two steps. First, extra structures are added on neural networks.
The resulting surrogate is a sum of two networks where one is the boundary network that is designed for fitting boundary data and the other is for fitting residual data.  
Importantly, the boundary network is pre-trained (or trained first) to fit the boundary data.
The rest of learning is done on the function that (approximately) satisfy the boundary conditions. 

PINNs \cite{Raissi_19_PINNs, Lagaris_00_ANN-Irregular} do require neither extra structures on the solution surrogate nor pre-training with the boundary data.
However, it has been empirically reported \cite{Lagaris_00_ANN-Irregular, Berg_18_Unified, Wang_20_GradPathPINN} that 
properly chosen boundary weights ($\lambda_{b}$) could accelerate the overall training and result in a better performance.
It was mentioned in \cite{Lagaris_00_ANN-Irregular} that the boundary weight takes a large positive value to accurately satisfy the boundary conditions.

On the top of these existing works, Theorem~\ref{cor:NN-elliptic} theoretically sheds light on the importance of learning the boundary conditions. 

Finally, we show that when  the solution to the PDE \eqref{def:elliptic-PDE}
is exactly represented by 
a neural network with a fixed architecture, 
one can further improve the mode of convergence.
\begin{corollary} \label{cor:NN-elliptic}
	Under the same conditions of Lemma~\ref{lem:NN} and Theorem~\ref{thm:main-elliptic}, 
	suppose there exists a class $\mathcal{H}_{\vec{\bm{n}}}^{\text{NN}}$
	of neural networks 
	with a fixed architecture
	\eqref{def:NN-class} 
	such that $u^* \in \mathcal{H}_{\vec{\bm{n}}}^{\text{NN}}$
	where $u^*$ is the solution to the PDE \eqref{def:elliptic-PDE}.
	Let $\mathcal{H}_{m_r} = \mathcal{H}_{\vec{\bm{n}}}^{\text{NN}}$ for all $m_r$.
	Then, with probability 1 over iid samples, 
	the sequence of minimizers stated in Theorem~\ref{thm:main-elliptic}
	converges to $u^*$ in $C^k(U)$,
	where $k$ is stated in Lemma~\ref{lem:NN}.
\end{corollary}
\begin{proof}
	By Theorem~\ref{thm:main-elliptic}, we have $\|h_{m_r} - u^*\|_{C^{0}(U)} \to 0$ as $m_r \to \infty$.
	Since $u^* \in \mathcal{H}_{\vec{\bm{n}}}^{\text{NN}}$,
	$u^* \in C^k(\overline{U})$ as well.
	It then follows from Lemma~\ref{lem:NN}
	that $\|h_{m_r} - u^*\|_{C^{k}(U)} \to 0$ as $m_r \to \infty$.
\end{proof}

%

\subsection{Linear Parabolic PDEs}
In this section, we consider the second-order linear parabolic equations.
Let $U$ be a bounded domain in $\mathbb{R}^{d}$
and let $U_T = U \times (0,T]$ for some fixed time $T > 0$.
Let us denote the parabolic boundary as $\Gamma_T = \overline{U}_T - U_T$.
Also, let 
$(\x,t)=(x_1,\cdots,x_d,t)$ be a point in $\mathbb{R}^{d+1}$.

Let us consider the initial/boundary-value problem:
\begin{equation} \label{def:Parabolic}
\begin{split}
\begin{cases}
-u_t+\mathcal{L}[u] = f,  & \text{in } U_T \\
u = \varphi, & \text{in } \partial U \times [0,T] \\
u = g, & \text{in } \overline{U} \times \{t=0\},
\end{cases}
\end{split}
\end{equation}
where 
$f:U_T \mapsto \mathbb{R}$, $g:\overline{U} \mapsto \mathbb{R}$, 
$\varphi:\partial U \times [0,T] \mapsto \mathbb{R}$,
and
\begin{align*}
\mathcal{L}[u] = \sum_{i,j=1}^{d} D_i\left(a^{ij}(\x,t) D_ju + b^i(\x,t)u\right) + \sum_{i=1}^{d} c^i(\x,t) D_iu + d(\x,t)u.
\end{align*} 

In order to ensure the existence, regularity and uniqueness of the PDE \eqref{def:Parabolic},
the following assumptions are made.
Again, we follow the Schauder approach \cite{Friedman_08_ParabolicPDEs}.
\begin{assumption} \label{assumption:parabolic}
	Let $\lambda(\x,t)$ be the minimum eigenvalues of $[a^{ij}(\x,t)]$
	and $\alpha \in (0,1)$.
	Suppose $a^{ij}, b^i$ are differentiable
	and let $\tilde{b}^i(\x,t) = \sum_{j=1}^d D_ja^{ij}(\x,t) + b^i(\x,t) + c^i(\x,t)$
	and $\tilde{c}(\x,t)  = \sum_{i=1}^{d} D_ib^i(\x,t)+d(\x,t)$.
	Let $\Omega = U\times (0,T)$.
	\begin{enumerate}
		\item For some constant $\lambda_0$, $\lambda(\x,t) \ge \lambda_0 > 0$
		for all $(\x,t) \in \Omega$.
		\item $a^{ij}, \tilde{b}^i, \tilde{c}^i$ are H$\ddot{o}$lder continuous (exponent $\alpha$) in $\Omega$,
		and $|a^{ij}|_{\alpha}, |\text{d}\tilde{b}^i|_{\alpha},
		|\text{d}^2\tilde{c}|_{\alpha}$
		are uniformly bounded.
		\item $f$ is H$\ddot{o}$lder continuous (exponent $\alpha$) in $U_T$
		and $|\text{d}^2f|_{\alpha} < \infty$.
		$g$ and $\varphi$ are H$\ddot{o}$lder continuous (exponent $\alpha$) on 
		$\overline{U} \times \{t=0\}$, and $\partial{U} \times [0,T]$, respectively, and $g(\x) = \varphi(\x,0)$ on $\partial U$.
		\item There exists $\theta > 0$ such that $\theta^2 a^{11}(\x,t) + \theta\tilde{b}^1(\x,t) \ge 1$ in $\Omega$.
		\item For $x' \in \partial U$, there exists a closed ball $\bar{B}$ in $\mathbb{R}^d$ such that $\bar{B}\cap \overline{U} = \{x'\}$.
		\item There are constants $\Lambda, v \ge 0$ such that for all $(\x, t) \in \Omega$
		\begin{equation*}
		    \sum_{i,j} |a^{ij}(\x,t)|^2 \le \Lambda, \qquad \lambda_0^{-2}(\|\bm{b}(\x,t)\|^2 + \|\bm{c}(\x,t)\|^2) + \lambda_0^{-1}|d(\x,t)| \le v^2,
		\end{equation*}
		where $\bm{b}(\x,t) = [b^1(\x,t),\cdots,b^d(\x,t)]^T$
		and $\bm{c}(\x,t) = [c^1(\x,t),\cdots, c^d(\x,t)]^T$.
	\end{enumerate}
\end{assumption}

By adopting the notation of Theorem~\ref{thm:main-elliptic}, 
we now present the convergence theorem for the linear parabolic PDEs.
\begin{theorem} \label{thm:main-parabolic}
	Suppose Assumptions~\ref{assumption:data-dist}, ~\ref{assumption:convergence}
	and ~\ref{assumption:parabolic} hold
	and $\tilde{c}(\x,t)  = \sum_{i=1}^{d} D_ib^i(\x,t)+d(\x,t) \le 0$.
	Let $m_r$ and $m_{b}$ be the number of iid samples from 
	$\mu_r$ and $\mu_{b}$, respectively,
	and $m_r = \mathcal{O}(m_{b}^{\frac{d+1}{d}})$.
	Let $\bm{\lambda}_{m_r}^R$ be a vector satisfying 
	\eqref{def:lambda-condition}.
	Let $h_{m_r} \in \mathcal{H}_{m_r}$ be a minimizer of the H\"{o}lder regularized loss $\text{Loss}_{m_r}(\cdot;\bm{\lambda}, \bm{\lambda}_{m_r}^R)$ \eqref{def:Holder-Reg-Loss}.
	Then the following holds:
	\begin{enumerate}
		\item For each PDE data set $\{f, g, \varphi\}$
		satisfying Assumption~\ref{assumption:parabolic}, 
		there exists a unique classical solution $u^*$ to the PDE \eqref{def:Parabolic}.
		\item With probability 1 over iid samples, 
		\begin{equation*}
		    \lim_{m_r \to \infty} h_{m_r} = u^*, \qquad \text{in } C^0(\Omega).
		\end{equation*}
	\end{enumerate}
\end{theorem}
\begin{proof}
	The proof can be found in Appendix~\ref{app:thm:main-parabolic}.
\end{proof}

{
Theorem~\ref{thm:main-parabolic} shows that neural networks that minimize the H\"{o}lder regularized empirical losses \eqref{def:Holder-Reg-Loss}
converge to the unique classical solution to the PDE \eqref{def:Parabolic}.
This answers the question
we posed in the introduction
for linear second-order parabolic PDEs.
}


The mode of convergence can be $L^2(0,T;H^1(U))$
if each minimizer satisfies
both the initial and the boundary conditions.
The Bochner space  $L^2(0,T;H^1(U))$ \cite{yosida1988functional}
consists
of all stronlgy measurable functions 
$u:[0,T] \mapsto H^1(U)$ with 
$\|u\|_{L^2(0,T;H^1(U))}^2 := \int_0^T \|u(t)\|_{H^1(U)}^2 < \infty$.
\begin{theorem} \label{thm:parabolic-H1}
	Under the same conditions of Theorem~\ref{thm:main-parabolic}, 
	suppose the sequence of minimizers of \eqref{def:problem}
	satisfies both the boundary and the initial conditions.
	Then, with probability 1 over iid samples, 
	the sequence of minimizers stated in Theorem~\ref{thm:main-parabolic} 
	strongly converges to the solution $u^*$ to \eqref{def:Parabolic}
	in $L^2(0,T;H^1(U))$.
\end{theorem}
\begin{proof}
	The proof can be found in Appendix~\ref{app:thm:parabolic-H1}.
\end{proof}
Again, Theorem~\ref{thm:parabolic-H1} shows the importance of learning the initial and boundary conditions
for the better convergence.

Similarly to Corollary~\ref{cor:NN-elliptic},
when the solution to the PDE \eqref{def:Parabolic} 
can be exactly represented by a class of neural networks with a fixed architecture, one could further improve the mode of convergence.
\begin{corollary}
	Under the same conditions of Lemma~\ref{lem:NN} and Theorem~\ref{thm:main-parabolic}, 
	suppose there exists a class $\mathcal{H}_{\vec{\bm{n}}}^{\text{NN}}$ of neural networks
	with a fixed architecture
	\eqref{def:NN-class} 
	such that 
	$u^* \in \mathcal{H}_{\vec{\bm{n}}}^{\text{NN}}$
	where $u^*$ is the solution to the PDE \eqref{def:Parabolic}.
	Let $\mathcal{H}_{m_r} = \mathcal{H}_{\vec{\bm{n}}}^{\text{NN}}$ for all $m_r$.
	Then, with probability 1 over iid samples, 
	the sequence of minimizers stated in Theorem~\ref{thm:main-parabolic} 
	converges to $u^*$ in $C^k(U_T)$,
	where $k$ is stated in Lemma~\ref{lem:NN}.
\end{corollary}
\begin{proof}
	Since the proof is similar to the proof of Corollary~\ref{cor:NN-elliptic},
	we omitted it.
\end{proof}

%% file: examples.tex
\section{Computational Examples}
\label{sec:example}

In this section, we provide computational examples to
illustrate our theoretical findings.
Mainly, we show 
the $L^2$ and $H^1$ convergence of the trained neural networks as the number of training data grows.
We confine ourselves to one or two-dimensional problems.

As an effort to find a neural network that minimizes the loss,
we employ a combination of two optimization methods; \texttt{Adam} \cite{Kingma_14_Adam}
and \texttt{L-BFGS} \cite{Liu_89_LBFGS}.
We consecutively apply these in the order of 
\texttt{Adam} and \texttt{L-BFGS}.
\texttt{Adam} is employed with its default hyper-parameter setting.
This combination has been used in other works (e.g. \cite{Lu_19_Deepxde}).
We employ the Xavier initialization \cite{Glorot2010understanding}
to initialize network parameters.

\subsection{Poisson's equation}
Let us consider the one-dimensional Poisson's equation:
\begin{equation*}
- u_{xx}(x) = f(x) \quad \text{in } U=(-1,1), \qquad u(-1) = g, u(1) = h.
\end{equation*}
The force term is simply obtained by substituting exact solution in the equation.
The training data are set to the equidistant points on $[-1,1]$. 
The prediction errors are computed based on the fixed grid of either 10,000 or 50,000 equidistant points on $[-1,1]$.
The discrete $L^2$ and $H^1$ predictive errors are reported. 
We show the training results by the original PINN loss \eqref{def:PINN-loss} and the Lipschitz regularized loss \eqref{def:LIPR-loss}.
We refer to the results by the original PINN loss as `PINN'
and to the results by the LIPschitz Regularized loss as `LIPR'.
Specifically, the results by `PINN'
and the results by `LIPR' are obtained by minimizing 
\begin{equation} \label{Poi-loss}
\begin{split}
\text{Loss}^{\text{PINN}}_{m_r}(h) &= \frac{1}{m_r}\sum_{j=1}^{m_r} \|\mathcal{L}[h](\x_r^j)- f(\x_r^j)\|^2
+ 
\frac{\lambda_b}{2}(|h(-1)|^2 +|h(1)|^2), \\
\text{Loss}^{\text{LIPR}}_{m_r}(h) &= \text{Loss}^{\text{PINN}}_{m_r}(h) + \lambda_r^R\max_{i} |\nabla \mathcal{L}[h](\x_r^i)|^2,
\end{split}
\end{equation}
respectively. 
Here $m_r$ is the number of residual training data points.

\begin{figure}[htbp]
	\centerline{
		\includegraphics[width=8cm]{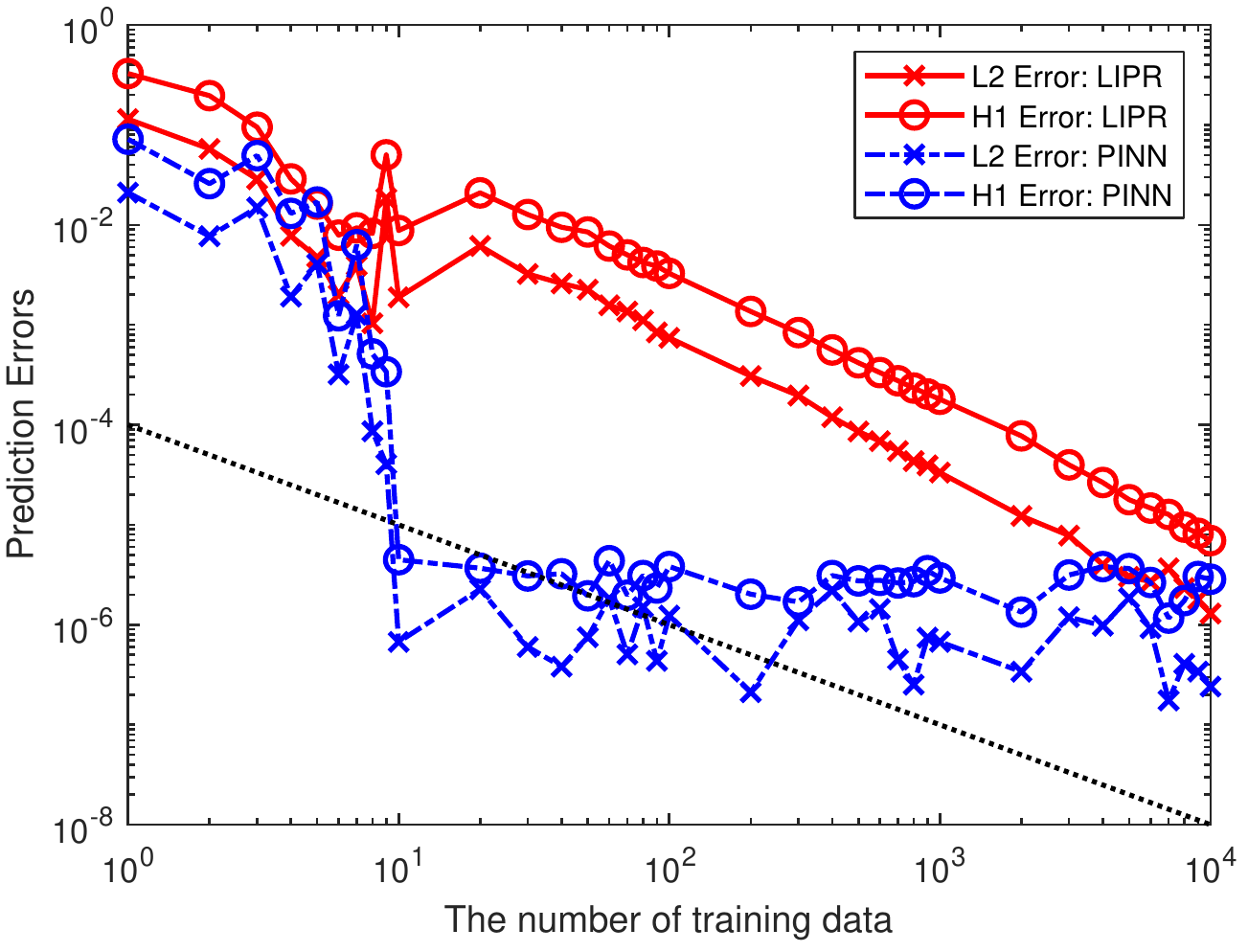}
	}
	\caption{The $L^2$ and $H^1$ convergence for the 1D Poisson equation whose exact solution is $u^*(x)=\tanh(x)$
	with respect to the number of training data points.
	The residual neural networks of depth 2 and width 50 are employed. 
	The `PINN' results are shown as dash-dot lines
	and the `LIPR' results are shown as solid lines.
	The dotted line is a reference line indicating the $\mathcal{O}(m_r^{-1})$-rate of convergence.
	}
	\label{fig:Poi1D-tanh}
\end{figure}

We first consider the case where the exact solution is given by 
$u^*(x) = \tanh(x)$. 
For this task, we employ the $2$-hidden layer residual neural network $u_{\text{NN}}$ having 50 neurons at each layer with the $\tanh$-activation function. 
That is,
$$
u_{\text{NN}}(x) = W^3\tanh(W^2(\bm{1}x+\tanh(W^1x + b^1)) + b^2 ) + b^3,
$$
where $W^j \in \mathbb{R}^{n_j\times n_{j-1}}$, $b^j \in \mathbb{R}^{n_j}$, $\bm{1} \in \mathbb{R}^{n_1}$ whose entries are all 1, and 
$\bm{\vec{n}} = (n_0=1, n_1=50, n_2=50, n_3 = 1)$.
We remark that in this case, the solution can be exactly represented by a neural network.
In Figure~\ref{fig:Poi1D-tanh}, 
the predictive errors are plotted with respect to the number of training data points from 1 to 10,000. 
We set $\lambda_b = 1$ and $\lambda_r^R = m_r^{-1.5}$ as suggested by 
Theorem~\ref{thm:gen}.
The maximum number of epochs of \texttt{Adam} is set to 25,000
and the full-batch is used. 
As expected by Corollary~\ref{cor:NN-elliptic}, 
the results by `LIPR' exhibit both $L^2$- and $H^1$-convergence of the errors.
We see that the rate of convergence approximately follows $\mathcal{O}(m_r^{-1})$,
which matches the rate of convergence of the expected PINN loss (Theorem~\ref{thm:conv-loss}).  
For reference, we plot a dotted-line showing the $\mathcal{O}(m_r^{-1})$ rate of convergence with respect to the number of training data.
The results by `PINN' converges very fast with respect to the number of training data and the errors are saturated at around 10 data points. 
Without any regularization terms, the `PINN' somehow finds an accurate approximation to $u^*$.
We observe that the predictive errors by `PINN' stay at the level of $10^{-6}$ at all time. 

\begin{figure}[htbp]
	\centerline{
		\includegraphics[width=8cm]{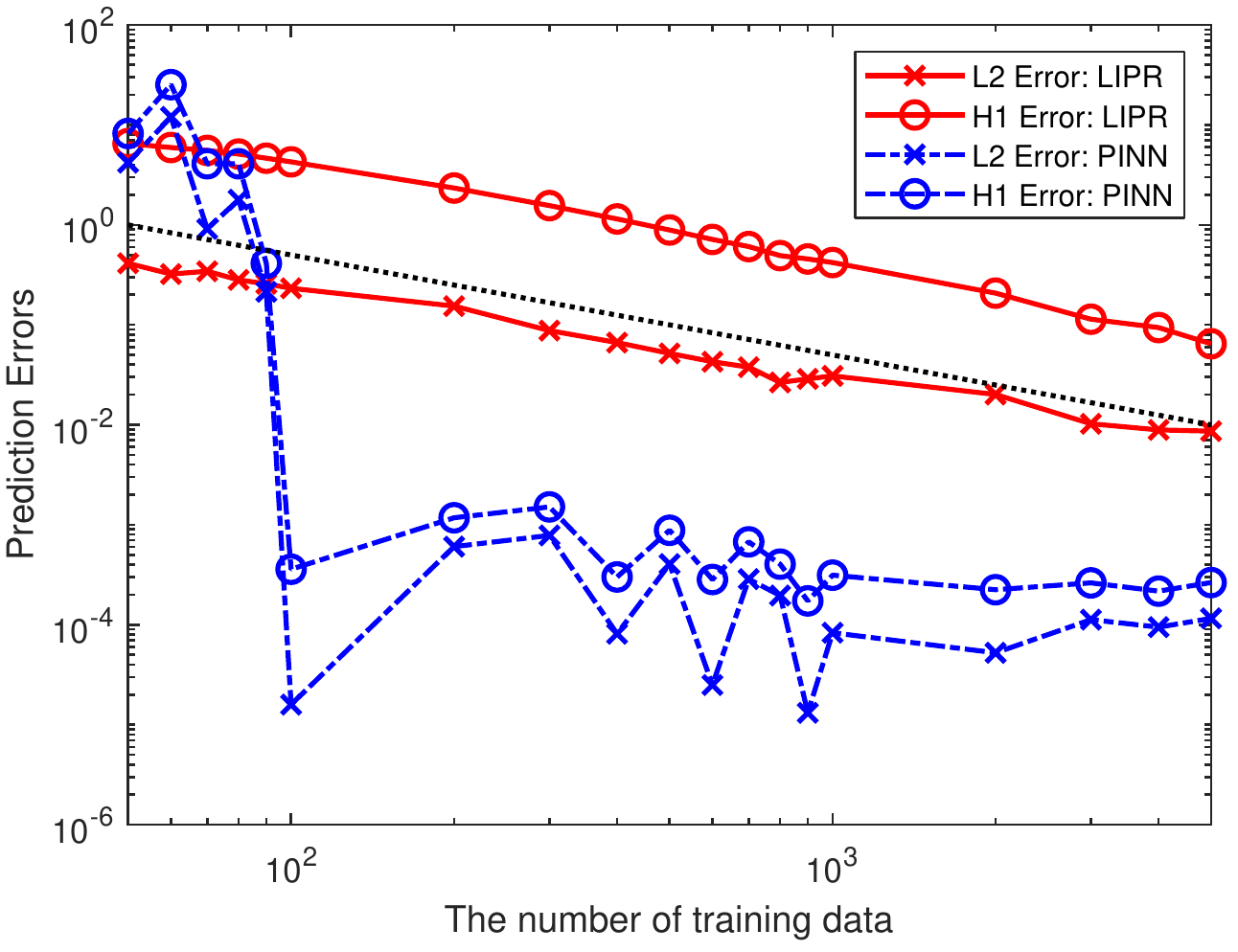}
	}
	\caption{The $L^2$ and $H^1$ convergence for the 1D Poisson equation whose exact solution is $u^*(x) = (1-x^2)\sin(6\pi x)$
	with respect to the number of training data points.
	The neural networks \eqref{def:bdry-match-NNs} that automatically satisfy the boundary conditions are employed.
	The `PINN' results are shown as dash-dot lines
	and the `LIPR' results are shown as solid lines.
	The dotted line is a reference line indicating the $\mathcal{O}(m_r^{-1})$-rate of convergence.
	}
	\label{fig:Poi1D-k6-H1}
\end{figure}

Next, we consider the case where the exact solution is given by 
$u^{*}(x) = (1-x^2)\sin(6\pi x)$.
For the training, we choose the function class that satisfies the boundary conditions automatically by adopting the following construction: 
\begin{equation} \label{def:bdry-match-NNs}
\hat{u}_{\text{NN}}(x) = (1-x^2)u_{\text{NN}}(x).
\end{equation}
Again, $u_{\text{NN}}$ is the $2$-hidden layer residual neural network having 50 neurons at each layer with the $\tanh$-activation function. 
Since the boundary conditions are automatically satisfied, the choice of $\lambda_b$ does not affect the training results.
We set $\lambda_r^R = {m_r}^{-1.5}$, the maximum number of epochs of \texttt{Adam} to either 5,000 or 25,000 and use the mini-batch of size 100. 
Figure~\ref{fig:Poi1D-k6-H1} shows the $L^2$- and $H^1$-predictive errors with respect to the number of training data points from 50 to 5,000. 
As expected by Theorem~\ref{thm:main-elliptic} and ~\ref{thm:elliptic-H1}, 
we see both the $L^2$- and $H^1$-convergence of the errors by `LIPR'.
Furthermore, we observe that the convergence rate approximately follows $m_r^{-1}$.
As a reference, the $\mathcal{O}(m_r^{-1})$ rate of convergence is plotted as a dotted line. We see that the slopes of the predictive errors by `LIPR' are well matched to
the slope of the dotted line.
The predictive errors by `PINN', again, decay very fast 
and are much smaller than those by 'LIPR' at all training data points but the first 100 points.
In this case, the predictive errors by `PINN' saturate at the level of $10^{-4}$.
This again indicates that PINNs find an approximation that also generalizes well without having explicit regularization terms. 

Although our experiments may not find global minima, 
our results clearly demonstrate both the $L^2$- and $H^1$-convergence of the errors
by the Lipschitz regulared empirical loss \eqref{def:LIPR-loss}.

\subsection{Heat Equation}
Let us consider the 1D heat equation:
\begin{equation*}
-u_t + \nu u_{xx} = f \quad \text{in } (x,t) \in U_T = (-1,1)\times (0,T],
\end{equation*}
with the Dirichlet initial/boundary condition.
There are two boundary domains ($\Gamma_1, \Gamma_2$) and one initial domain ($\Gamma_3$): 
$\Gamma_1 = \{-1\} \times [0,T]$,
$\Gamma_2 = \{1\} \times [0,T]$, and $\Gamma_3 = (-1,1) \times \{0\}$.
We consider the case where the exact solution is given by
$u^*(x,t) = \sin(\pi x)e^{-t}$
and $T = 1$.
Let $m_{b,j}$ be the number of training points on $\Gamma_j$
and $m_r$ be the number of residual points on $U_T$.
The training points are randomly uniformly drawn from its corresponding domains.
We set $m_{b,2} = m_{b,1}$, $m_{b,3} = 2m_{b,1}$ and $m_r = 2m_{b,1}m_{b,2}$.
This satisfies the condition of $m_r = \mathcal{O}(m_{b}^2)$, which is stated in Theorem~\ref{thm:main-parabolic}.
Again, we report two training results: the original PINN loss \eqref{def:PINN-loss} (`PINN')
and the Lipschitz regularized loss \eqref{def:LIPR-loss} (`LIPR').
For the `PINN' loss, we set all the weights to 1.
For the `LIPR' weights, we set $\bm{\lambda} = \bm{1}$,  
$\lambda_{r}^R = \frac{2}{m_r}$, $\lambda_{b,1}^R = \frac{1}{m_{b,1}\sqrt{m_r}}$,
$\lambda_{b,2}^R = \frac{1}{m_{b,2}\sqrt{m_r}}$, and $\lambda_{b,3}^R = \frac{1}{m_{b,3}\sqrt{m_r}}$.
The standard feed-forward $\tanh$-neural networks of depth 2 and width 50 are employed for the experiments.
We set the maximum number of epoch of \texttt{Adam} to be 10,000 and use 
the mini-batch training.

\begin{figure}[htbp]
	\centerline{
		\includegraphics[width=8cm]{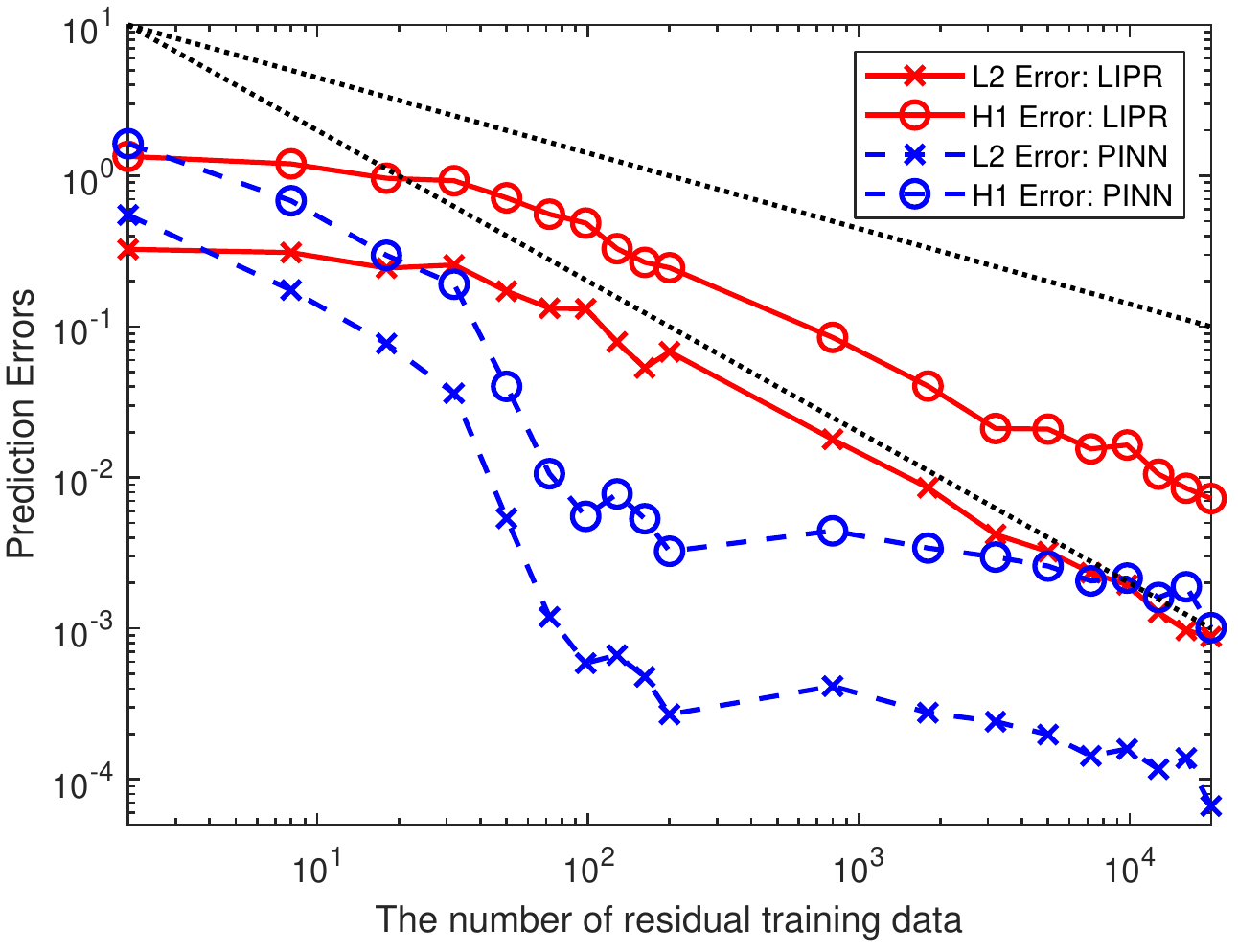}
	}
	\caption{
		The $L^2(0,1;L^2)$- and $L^2(0,1;H^1)$-convergence of the errors with respect to the number of training data point 
		for the 1D Heat equation whose exact solution is $u^*(x) = \sin(\pi x)e^{-t}$.
		The feed-forward neural networks of depth 2 and width 50 are employed.
		The `PINN' results are shown as dash-dot lines
		and the `LIPR' results are shown as solid lines.
		The dotted lines are reference lines indicating the $\mathcal{O}(m_r^{-\frac{1}{2}})$- and $\mathcal{O}(m_r^{-1})$-rate of convergence.
	}
	\label{fig:Heat1D-L2}
\end{figure}

In Figure~\ref{fig:Heat1D-L2}, we show the predictive $L^2(0,1;L^2)$- and $L^2(0,1;H^1)$-errors with respect to the increasing the number of residual training data points $m_r$.
The number of residual training data increases according to $m_r = 2m_{b,1}^2$ for $m_{b,1}=1, 2, \cdots, 10, 20, \cdots,100$. 
The predictive errors are computed on the tensor-grid constructed 
by using 400 and 200 equidistant points on $[-1,1]$ and $[0,1]$, respectively.
We report the average of three independent simulations. 
We see that the results by `LIPR' exhibit both the $L^2(0,1;L^2)$- and $L^2(0,1;H^1)$-convergence.
We also observe that the rate of the $L^2(0,1;L^2)$-convergence by `LIPR' approximately follows $\mathcal{O}(m_r^{-1})$.
For a reference, we plot a line showing the $\mathcal{O}(m_r^{-1})$-rate of convergence
as dotted line starting at $(200, 0.1)$.
It can be seen that the slope of the $L^2(0,1;L^2)$-errors by `LIPR' is well matched to those of the dotted line.
We see that the rate of the $L^2(0,1;H^1)$-convergence by `LIPR' lies between $\mathcal{O}(m_r^{-1})$
and $\mathcal{O}(m_r^{-\frac{1}{2}})$. 
Again, we plot a line showing the $\mathcal{O}(m_r^{-\frac{1}{2}})$-rate of convergence
as dotted line starting at $(200, 0.2)$.
It follows from Theorem~\ref{thm:conv-loss} that in this case, the expected PINN loss decays at the rate at least $\mathcal{O}(m_r^{-\frac{1}{2}})$.
We observe that this rate continues to hold for the predictive errors as demonstrated in Figure~\ref{fig:Heat1D-L2}.
The results by `PINN' converge faster and produce smaller predictive errors compared to those by `LIPR' at almost all training data points.
This again indicates that the empirical PINN loss is capable of producing a good approximation to the solution without having an explicit regularization. 

%
%
%
%

%% file: conclusion.tex
\section{Conclusion} \label{sec:conclusion}

{
In this work, we prove the consistency 
of physics informed neural networks (PINNs)
in the sample limit. 
It is equivalent to the convergence of generalization error in the context of machine learning.
Upon deriving an upper bound of the expected PINN loss, 
we obtain a H\"{o}lder regularized empirical loss. 
Under some assumptions, we show that the expected PINN loss at minimizers of the H\"{o}lder regularized loss converges to zero.
By considering two classes of PDEs --linear second-order elliptic and parabolic that admit classical solutions--
we show that with probability 1 over iid training samples, 
a sequence of neural networks that minimize the H\"{o}lder regularized losses converges to the solution to the PDE uniformly.
Furthermore, we show that if each minimizer exactly satisfies the boundary conditions, the mode of convergence becomes $H^1$. 
We provided a set of conditions 
for neural networks 
on which the convergence can be guaranteed.
}





%% file: appendix.tex
\section{Proof of Lemma~\ref{lem:NN}} \label{app:lem:NN}
\begin{proof}
	For each $h_j$, let $\bm{\theta}^j$ be its associated weights and biases.
	Since $\{\bm{\theta}^j\}$ is uniformly bounded,
	there exists a convergent subsequence, say, $\{\bm{\theta}^{j_i}\}$.
	Let $\{h_{j_i}\}$ be its corresponding neural network sequence.
	Let $\bm{\theta}^*$ be the limit of $\{\bm{\theta}^{j_i}\}$ and 
	$h^* \in \mathcal{H}_{\vec{\bm{n}}}^{NN}$ be the corresponding limit network.
	Note that since $\{\bm{\theta}_j\}$ is uniformly bounded, $U$ is bounded,
	and $\sigma^{(s)}$ is bounded for all $s \in \{0,\cdots,k\}$, 
	$D^{\bm{s}} h_{j_i}$ is bounded continuous in $\bar{U}$
	for all $\bm{s}$ with $0 \le |\bm{s}|\le k$.
	
	\begin{claim} \label{claim-1}
		Let $\{h^L(\x;\bm{\theta}^{j_i})\}$ be a sequence of $L$-layer neural networks
		having the same architecture $\vec{\bm{n}}$
		whose activation function $\sigma$ is of $C^k(\mathbb{R})$
		and  
		$\frac{d^s \sigma(x)}{dx^s}$ is bounded and Lipschitz continuous for $s=0,\cdots,k$.
		Let $\bm{s} = (s_1,\cdots, s_{d})$ be a multi-index with
		$|\bm{s}|=\sum_{i=1}^d s_i$.
		Then, for each $\bm{s}$ with $0 \le |\bm{s}|\le k$, 
		the sequence $\{D^s h^L(\x;\bm{\theta}^{j_i})\}$ is uniformly convergent to
		$D^s h^L(\x;\bm{\theta}^{*})$ on $\bar{U}$.
	\end{claim}
	\begin{proof}[Proof of Claim~\ref{claim-1}]
		Let $\vec{\bm{n}} = (n_0,\cdots, n_L)$ be the network architecture.
		Let us recall that given $\vec{\bm{n}}$, 
		$$
		h^l(\x) = \bm{W}^l\sigma(h^{l-1}(\x)) + \bm{b}^l, \qquad 1 < l \le L,
		$$
		and $h^1(\x) = \bm{W}^1\x + \bm{b}^1$,
		where $\bm{W}^l \in \mathbb{R}^{n_l\times n_{l-1}}$ and $\bm{b}^l \in \mathbb{R}^{n_l}$.
		For each $l$, let $h^l_k$ be the $k$-th output of $h^l$,
		where $1\le k \le n_l$.
		Let $\bm{\theta}_l = \{\bm{W}^j, \bm{b}^j\}_{j=1}^l$.
		With a slight abuse of notation, we regard $\bm{\theta}_l$ as a column vector.
		Also, let $\bm{\theta'}_l = \{\bm{W'}^j, \bm{b'}^j\}_{j=1}^l$.
		
		We want to prove the following statement for all $l =1,\cdots, L$:
		\begin{align*}
		D^{\bm{s}} h^l(\x;\bm{\theta}^{j_i})
		\text{ converges to }
		D^{\bm{s}} h^l(\x;\bm{\theta}^*)
		\text{ in } C^0(U),
		\forall \bm{s} \text{ with } 0 \le |\bm{s}| \le k.
		\end{align*}
		We prove it by applying induction on $l$.

		When $l=1$, it is clear that for $r=1,\cdots,n_1$,
		we have
		\begin{align*}
		|h_r^1(\x;\bm{\theta}_1) - h_r^1(\x;\bm{\theta}_1')|^2
		&\le 2\left(|\bm{b}^1_r - \bm{b'}^1_r|^2 + \|\bm{W}^1_r - \bm{W'}^1_r\|^2\|x\|^2\right)
		\\
		&\le C\|\bm{\theta}_1 - \bm{\theta'}_1\|^2,
		\end{align*}
		where $C = 2\max\{1, \sup_{U} \|x\|^2\}$, and
		for $\bm{s}=(s_1,\cdots,s_d)$ with $1\le |\bm{s}| \le k$, 
		\begin{align*}
		|D^{\bm{s}}h_r^1(\x;\bm{\theta}) - D^{\bm{s}}h_r^1(\x;\bm{\theta}')|
		&= \begin{cases}
		0 & \text{if } \max s_i > 1, \\
		|\bm{W}^1_{ri} - \bm{W'}^1_{ri}| & \text{if } s_i = 1.
		\end{cases}
		\end{align*}
		Hence, for $\bm{s}=(s_1,\cdots,s_d)$ with $0 \le |\bm{s}| \le k$,
		$D^{\bm{s}} h^1(\x;\bm{\theta}_1^{j_i})$ converges to
		$D^{\bm{s}} h^1(\x;\bm{\theta}_1^*)$ in $C^0(U)$.

		Suppose the statement is true for $l-1$ where $l \ge 2$ 
		and we want to show the case for $l$.
		For $\bm{s}$ with $0 \le |\bm{s}| \le k$,
		it then can be checked that for $1\le r \le n_l$,
		\begin{align*}
		D^{\bm{s}} h^l_r(\x;\bm{\theta}_l) = \frac{\partial^{|\bm{s}|} h_r^l(\x;\bm{\theta}_l)}{\partial x_1^{s_1}\cdots \partial x_d^{s_d}} =
		\sum_{i=1}^{n_l} \bm{W}^l_{ri}
		\sum_{t=0}^{|\bm{s}|}
		Q_{i,t}^{l-1}(\bm{\theta}_{l-1},|\bm{s}|,\x) \sigma^{(t)}(h^{l-1}_i(\x;\bm{\theta}_{l-1})),
		\end{align*}
		where $\sigma^{(t)}$ is the $t$-th derivative of $\sigma$ and
		$\{Q_{i,t}^l(\bm{\theta}_l,|\bm{s}|,\bm{x})\}_{t=1}^{|\bm{s}|}$ is recursively defined as follows:
		Let 
		$\bar{x}_{j} = x_i$ if $\sum_{l=1}^{i-1} s_l < j \le \sum_{l=1}^i s_l$.
		Then, for $s \ge 2$, 
		\begin{align*}
		&Q_{i,t}^{l-1}(\bm{\theta}_{l-1},s,\bm{x}) \\
		&=
		\begin{cases}
		Q_{i,s-1}^{l-1}(\bm{\theta}_{l-1},s-1,\bm{x})\frac{\partial h^{l-1}_i(\x;\bm{\theta}_{l-1})}{\partial \bar{x}_s}
		& \text{if } t = s, \\
		\frac{\partial}{\partial \bar{x}_s}Q_{i,t}^{l-1}(\bm{\theta}_{l-1},s-1,\bm{x})
		+ Q_{i,t-1}^{l-1}(\bm{\theta}_{l-1},s-1,\bm{x})\frac{\partial h^{l-1}_i(\x;\bm{\theta}_{l-1})}{\partial \bar{x}_s}
		& \text{if } 1 < t < s, \\
		\frac{\partial}{\partial \bar{x}_s}Q_{i,1}^{l-1}(\bm{\theta}_{l-1},s-1,\bm{x})
		&\text{if } t= 1, \\
		0 & \text{if } t =0,
		\end{cases}
		\end{align*}
		with $Q_{i,1}^{l-1}(\bm{\theta}_{l-1},1,\bm{x}) =\frac{\partial h^{l-1}_i(\x;\bm{\theta}_1)}{\partial \bar{x}_1}$,
		$Q_{i,0}^{l-1}(\bm{\theta}_{l-1},1,\bm{x}) = 0$,
		and $Q_{i,0}^{l-1}(\bm{\theta}_{l-1},0,\bm{x}) = 1$.
		Thus,
		\begin{align*}
		&|D^{\bm{s}}h_r^l(\x;\bm{\theta}_{l-1}) - D^{\bm{s}}h_r^l(\x;\bm{\theta}_{l-1}')|
		\\
		&\le 
		\sum_{i=1}^{n_l} 
		\sum_{t=0}^{|\bm{s}|}
		|\bm{W}^l_{ri}Q_{i,t}^{l-1}(\bm{\theta}_{l-1},|\bm{s}|,\x) - \bm{W'}^l_{ri}Q_{i,t}^{l-1}(\bm{\theta'}_{l-1},|\bm{s}|,\x)|
		|\sigma^{(t)}(h^{l-1}_i(\x;\bm{\theta}_{l-1}))|
		\\
		&\quad+
		\sum_{i=1}^{n_l} 
		\sum_{t=0}^{|\bm{s}|} 
		|\bm{W'}^l_{ri}Q_{i,t}^{l-1}(\bm{\theta'}_{l-1},|\bm{s}|,\x)|
		|\sigma^{(t)}(h^{l-1}_i(\x;\bm{\theta}_{l-1}))-\sigma^{(t)}(h^{l-1}_i(\x;\bm{\theta'}_{l-1}))|.
		\end{align*}
		Since $\sigma^{(s)}$ is Lipchitz continuous and bounded for all $0 \le s \le k$,
		there exists two constants $M$ and $L$ such that  
		\begin{align*}
		|\sigma^{(t)}(h^{l-1}_i(\x;\bm{\theta}_{l-1}))-\sigma^{(t)}(h^{l-1}_i(\x;\bm{\theta'}_{l-1}))|
		&\le L|h^{l-1}_i(\x;\bm{\theta}_{l-1}) - h^{l-1}_i(\x;\bm{\theta'}_{l-1})|,
		\\
		|\sigma^{(t)}(h^{l-1}_i(\x;\bm{\theta}_{l-1}))| &\le M.
		\end{align*}
		Let $G^l(\x;\bm{\theta}_l) = \bm{W}^l_{ri}Q_{i,t}^{l-1}(\bm{\theta}_{l-1},|\bm{s}|,\x)$.
		By the induction hypothesis, $D^{\bm{s}} h^{l-1}(\x;\bm{\theta}_{l-1}^{j_i})$ converges to
		$D^{\bm{s}} h^{l-1}(\x;\bm{\theta}_{l-1}^*)$ in $C^0(U)$
		for any $\bm{s}$ with $0\le |\bm{s}| \le k$.
		Recall that the product of two uniformly convergent sequences of bounded continuous functions is also uniformly convergent to the product of limits.
		Thus, 
		$G^l(\x;\bm{\theta}_l^{j_i})$ uniformly converges to
		$G^l(\x;\bm{\theta}_l^*)$ in $U$.
		Therefore, by combining all the above,
		we conclude that 
		$D^{\bm{s}} h^l(\x;\bm{\theta}_{l}^{j_i})$ converges to
		$D^{\bm{s}} h^l(\x;\bm{\theta}_l^*)$ in $C^0(U)$.
		By induction, the proof is completed.
	\end{proof}
	
	By Claim~\ref{claim-1}, we conclude that $\lim_{i \to \infty} \|h_{j_i} - h^*\|_{C^k(U)} = 0$.
	Since $\lim_{j \to \infty} \|h_j - u\|_{C^0(U)} = 0$,
	we have $h^* = u$. 
	Hence, $\lim_{i \to \infty} \|h_{j_i} - u\|_{C^k(U)} = 0$.
	Since the subsequence was chosen arbitrarily, the proof is completed.
\end{proof}

\section{Proof of Theorem~\ref{thm:gen}} \label{app:thm:gen}
The proof consists of two lemmas.
\begin{lemma} \label{app:lemma-lip}
	Let $\mathcal{T}_r = \{\x_r^i\}_{i=1}^{m_r}$
	and $\mathcal{T}_{b} = \{\x_{b}^i\}_{i=1}^{m_{b}}$.
	Suppose $m_r$ and $m_{b}$ are large enough to satisfy the following:
	for any $\x_r \in U$ and $\x_{b} \in \Gamma$,
	there exists $\x_r' \in \mathcal{T}_r$ and $\x_{b}' \in \mathcal{T}_{b}$ such that 
	$
	\|\x_r - \x_r' \| \le \epsilon_r$ and 
	$\|\x_{b} - \x_{b}' \| \le \epsilon_{b}$.
	Then, we have
	\begin{equation} \label{app:lemma-pf-main-eqn}
	\begin{split}
    \text{Loss}^{\text{PINN}}(h;\bm{\lambda})
    &\le
	C_{\bm{m}} \cdot \left[\text{Loss}_{\bm{m}}^{\text{PINN}}(h;\bm{\lambda})
	+\frac{3\lambda_r\epsilon_r^{2\alpha}}{C_{\bm{m}}}\big[\mathcal{L}[h]\big]_{\alpha;U}^2
	+\frac{3\lambda_{b}\epsilon_{b}^{2\alpha}}{C_{\bm{m}}}\big[\mathcal{B}[h]\big]_{\alpha;\Gamma}^2\right]
	\\
	&\qquad
	+3\lambda_r\epsilon_r^{2\alpha}\big[f\big]_{\alpha;U}^2 +3\lambda_{b}\epsilon_{b}^{2\alpha}\big[g\big]_{\alpha;\Gamma}^2,
	\end{split}
	\end{equation}
	where $C_r$ and $C_{b}$ are from Assumption~\ref{assumption:data-dist}
	and
	$
	C_{\bm{m}} = 3\max\left\{C_rm_r\epsilon_r^{d},  C_{b}m_{b}\epsilon_{b}^{d-1} \right\}.
	$
\end{lemma}
\begin{proof}
	Let us first recall that we have
	$
	\|x + y + z\|^2 \le 3(\|x\|^2 + \|y\|^2 + \|z\|^2)
	$ for any three vectors $x, y, z$.
	For $\x_r, \x_r' \in U$, we have
	\begin{align*}
	&\|\mathcal{L}[h](\x_r) - f(\x_r)\|^2 
	\le 3\left(\|\mathcal{L}[h](\x_r) -\mathcal{L}[h](\x_r')\|^2 + \|\mathcal{L}[h](\x_r') -f(\x_r')\|^2 + \|f(\x_r') - f(\x_r)\|^2\right).
	\end{align*}
	Similarly, for $\x_{b}, \x_{b}' \in \Gamma$, we have
	\begin{align*}
	\|\mathcal{B}[h](\x_{b}) - g(\x_{b})\|^2
	&\le 3
	\|\mathcal{B}[h](\x_{b})- \mathcal{B}[h](\x_{b}')\|^2
	+3\|\mathcal{B}[h](\x_{b}') - g(\x_{b}')\|^2 
	 +3\|g(\x_{b}')- g(\x_{b})\|^2.
	\end{align*}
	By assumption, 
	$\forall \x_r \in U$ and $\forall \x_{b} \in \Gamma$,
	there exists $\x_r' \in \mathcal{T}_r^{m_r}$
	and $\x_{b}' \in \mathcal{T}_b^{m_{b}}$ such that 
	$\|\x_r - \x_r'\| \le \epsilon_r$ and $\|\x_{b} - \x_{b}'\| \le \epsilon_{b}$.
	Thus, we have
	\begin{align*}
	&\mathbf{L}(\x_r,\x_b;h,\bm{\lambda},\bm{0}) =
	\lambda_r\|\mathcal{L}[h](\x_r) - f(\x_r)\|^2 
	+ \lambda_{b}\|\mathcal{B}[h](\x_{b}) - g(\x_{b})\|^2
	\\
	&\le 3\lambda_r\|\mathcal{L}[h](\x_r) -\mathcal{L}[h](\x_r')\|^2 +
	3\lambda_r\|f(\x_r') - f(\x_r)\|^2
	+ 3\lambda_r\|\mathcal{L}[h](\x_r') - f(\x_r')\|^2
	\\
	&\quad
	+3\lambda_{b}\|\mathcal{B}[h](\x_{b})- \mathcal{B}[h](\x_{b}')\|^2
	+3\lambda_{b}\|g(\x_{b}')- g(\x_{b})\|^2
	+3\lambda_{b}\|\mathcal{B}[h](\x_{b}') - g(\x_{b}') \|^2
	\\
	&\le 3\textbf{L}(\x_r',\x_b';h,\bm{\lambda},\bm{0})
	+3\lambda_r\epsilon_r^{2\alpha}\big[\mathcal{L}[h]\big]_{\alpha;U}^2
	+3\lambda_{b}\epsilon_{b}^{2\alpha}\big[\mathcal{B}[h]\big]_{\alpha;\Gamma}^2 
	+3\lambda_r\epsilon_r^{2\alpha}\big[f\big]_{\alpha;U}^2 +3\lambda_{b}\epsilon_{b}^{2\alpha}\big[g\big]_{\alpha;\Gamma}^2.
	\end{align*}

	
	For $\x_r^i \in \mathcal{T}_r^{m_r}$,
	let $A_{\x_r^i}$ be the Voronoi cell assciated with $\x_r^i$,
	i.e., 
	$$
	A_{\x_r^i}= \{\x \in U | \|\x - \x_r^i\| = \min_{\x' \in \mathcal{T}_r^{m_r}} \|\x - \x'\| \},
	$$
	and let $\gamma^i_{r} = \mu_{r}(A_{\x^i_{r}})$.
	Similarly, for $\x^i_{b} \in \mathcal{T}_{b}^{m_{b}}$, 
	let 
	$$
	A_{\x^i_{b}} = \{\x \in \Gamma | \|\x - \x^i_{b}\| = \min_{\x' \in \mathcal{T}_{b}^{m_{b}}} \|\x - \x'\| \},
	$$
	and let $\gamma^i_{b} = \mu_{b}(A_{\x^i_{b}})$.
	Note that $\sum_{i=1}^{m_r} \gamma^i_f = 1$ and 
	$\sum_{i=1}^{m_{b}} \gamma^i_{b} = 1$.
	By taking the expectation with respect to $(\x_r, \x_b) \sim \mu = \mu_r\times \mu_{b}$,
	we have
	\begin{equation*} 
	\begin{split}
	&\mathbb{E}_{\mu}[\mathbf{L}(\x_r,\x_b;h,\bm{\lambda},\bm{0})]
	\le
	3\sum_{j=1}^{m_{b}}\sum_{i=1}^{m_r}
	\gamma_r^i\gamma_{b}^j
	\mathbf{L}(\x_r^i,\x_{b}^j;h,\bm{\lambda},\bm{0})
	\\
	&\qquad
	+3\lambda_r\epsilon_r^{2\alpha}\big[\mathcal{L}[h]\big]_{\alpha;U}^2
	+3\lambda_{b}\epsilon_{b}^{2\alpha}\big[\mathcal{B}[h]\big]_{\alpha;\Gamma}^2
	+3\lambda_r\epsilon_r^{2\alpha}\big[f\big]_{\alpha;U}^2 +3\lambda_{b}\epsilon_{b}^{2\alpha}\big[g\big]_{\alpha;\Gamma}^2.
	\end{split}
	\end{equation*}
	By letting $\gamma_r^{m_r,*} = \max_i \gamma_r^i$ and $\gamma_{b}^{m_{b},*} = \max_i \gamma^i_{b}$,
	we obtain
	\begin{align*}
	&\mathbb{E}_{\mu}[\mathbf{L}(\x_r,\x_b;h,\bm{\lambda},\bm{0})]
	\\
	&\le
	3m_r \gamma_r^{m_r,*}  \cdot \frac{\lambda_r}{m_r}\sum_{i=1}^{m_r} \|\mathcal{L}[h](\x_r^i) - f(\x_r^i)\|^2
	+
	3m_{b}\gamma_{b}^{m_{b},*} 
	\cdot \frac{\lambda_{b} }{m_{b}}\sum_{j=1}^{m_{b}} \|\mathcal{B}[h](\x_{b}^j) - g(\x_{b}^j)\|^2
	\\
	&\qquad 
	+3\lambda_r\epsilon_r^{2\alpha}\big[\mathcal{L}[h]\big]_{\alpha;U}^2
	+3\lambda_{b}\epsilon_{b}^{2\alpha}\big[\mathcal{B}[h]\big]_{\alpha;\Gamma}^2
	+3\lambda_r\epsilon_r^{2\alpha}\big[f\big]_{\alpha;U}^2 +3\lambda_{b}\epsilon_{b}^{2\alpha}\big[g\big]_{\alpha;\Gamma}^2.
	\end{align*}
	Note that $m_r\gamma_r^{m_r,*}, m_{b}\gamma_{b}^{m_{b},*} \ge 1$.
	Let $B_\epsilon(\x)$ be a closed ball centered at $\x$ with radius $\epsilon$.
	Let $P_r^* = \max_{\x \in U} \mu_r(B_{\epsilon_r}(\x) \cap U)$
	and $P_{b}^* = \max_{\x \in \Gamma} \mu_{b}(B_{\epsilon_{b}}(\x) \cap \Gamma)$.
	Since for any $\x_r \in U$ and $\x_{b} \in \Gamma$,
	there exists $\x_r' \in \mathcal{T}_r$ and $\x_{b}' \in \mathcal{T}_{b}$ such that 
	$\|\x_r - \x_r' \| \le \epsilon_r$, and
	$\|\x_{b} - \x_{b}' \| \le \epsilon_{b}$ 
	for each $i$, 
	there are closed balls $B_{\epsilon_r}$ and $B_{\epsilon_{b}}$
	that include $A_{\x_r^i}$ and $A_{\x_{b}^i}$, respectively.
	Thus, we have $\gamma_r^{m_r,*} \le P_r^*, 
	\gamma_{b}^{m_{b},*} \le P_{b}^*$.
	Moreover, it follows from  Assumption~\ref{assumption:data-dist} that 
	\begin{equation} \label{gamma-upper-bd}
	\gamma_r^{m_r,*} \le P_r^* \le C_r\epsilon_r^{d},
	\qquad
	\gamma_{b}^{m_{b},*} \le P_{b}^* \le C_{b}\epsilon_{b}^{d-1}.
	\end{equation}
	Therefore, we obtain
	\begin{equation*} 
	\begin{split}
	&\mathbb{E}_{\mu}[\mathbf{L}(\x_r,\x_b;h,\bm{\lambda},\bm{0})]
	\\ 
	&\le
	3C_rm_r\epsilon_r^{d}  \cdot \frac{\lambda_r}{m_r}\sum_{i=1}^{m_r} \|\mathcal{L}[h](\x_r^i) - f(\x_r^i)\|^2
	+3\lambda_r\epsilon_r^{2\alpha}\big[\mathcal{L}[h]\big]_{\alpha;U}^2
	+3\lambda_r\epsilon_r^{2\alpha}\big[f\big]_{\alpha;U}^2
	\\
	&\quad+
	3C_{b}m_{b}\epsilon_{b}^{d-1} 
	\cdot \frac{\lambda_{b} }{m_{b}}\sum_{j=1}^{m_{b}} \|\mathcal{B}[h](\x_{b}^j) - g(\x_{b}^j)\|^2
	+3\lambda_{b}\epsilon_{b}^{2\alpha}\big[\mathcal{B}[h]\big]_{\alpha;\Gamma}^2
	+3\lambda_{b}\epsilon_{b}^{2\alpha}\big[g\big]_{\alpha;\Gamma}^2
	\\
	&\le
	C_{\bm{m}} \cdot
	\text{Loss}_{\bm{m}}(h;\bm{\lambda},\bm{0})
	+3\lambda_r\epsilon_r^{2\alpha}\big[\mathcal{L}[h]\big]_{\alpha;U}^2
	+3\lambda_{b}\epsilon_{b}^{2\alpha}\big[\mathcal{B}[h]\big]_{\alpha;\Gamma}^2
	+3\lambda_r\epsilon_r^{2\alpha}\big[f\big]_{\alpha;U}^2 +3\lambda_{b}\epsilon_{b}^{2\alpha}\big[g\big]_{\alpha;\Gamma}^2,
	\end{split}
	\end{equation*}
	where
	$
	C_{\bm{m}} = 3\max\{C_rm_r\epsilon_r^{d},
	C_{b}m_{b}\epsilon_{b}^{d-1} \}.
	$
	Then, the proof is completed.
\end{proof}
\begin{lemma} \label{cor-sampling}
	Let $X$ be a compact subset in $\mathbb{R}^d$.
	Let $\mu$ be a probability measure supported on $X$.
	Let $\rho$ be the probability density of $\mu$ with respect to 
	$s$-dimensional Hausdorff measure on $X$
	such that $\inf_{X} \rho > 0$.
	Suppose that 
	for $\epsilon > 0$, 
	there exists a partition of $X$, $\{X_k^\epsilon\}_{k=1}^{K_{\epsilon}}$
	that depends on $\epsilon$ 
	such that for each $X_k^\epsilon$, 
	$c\epsilon^s\le \mu(X_k^\epsilon)$ where $c > 0$ depends only on $(\mu, X)$,
	and 
	there exists a cube $H_{\epsilon}(\textbf{z}_k)$ of side length $\epsilon$ centered at some $\textbf{z}_k$ in $X_k$
	such that $X_k \subset H_{\epsilon}(\textbf{z}_k)$.
	Then,
	with probability at least $1 - \sqrt{n}(1-1/\sqrt{n})^n$ over iid $n$ sample points $\{\x_i\}_{i=1}^n$ from $\mu$,
	for any $\x \in X$, there exists a point $\x_j$
	such that $\|\x - \x_j\| \le \sqrt{d}c^{-\frac{1}{s}}n^{-\frac{1}{2s}}$.
\end{lemma}
\begin{proof}
	From the conditions on $\mu$ and $X$, 
	$K:=K_\epsilon$ can be at most $\frac{1}{c\epsilon^{s}}$.
	Let $X_i := X_i^\epsilon$.
	Note that $\mu(X_i\cap X_j) = 0$ for $i\ne j$
	and $\mu(\cup_{j=1}^K X_j) = 1$.
	Let $p_i = \mu(X_i)$,
	the probability that a random sample from $\mu$ falls in $X_i$.
	Moreover, by the property of $\mu$, we have
	$p_i \ge c\epsilon^s$.
	
	For a positive integer $n$ satisfying $n \ge \frac{1}{c\epsilon^s} \ge K$,
	let $A_n$ be the event that for randomly drawn $n$ points, 
	each $X_i$ contains at least one point.
	Then,
	$$
	\mu(A_n) = \sum_{i_j \ge 1, |\bm{i}|=n} \binom{n}{i_1,\cdots,i_K}  p_1^{i_1}\cdots p_K^{i_K},
	$$
	where $\bm{i} = (i_1,\cdots,i_K)$ and $|\bm{i}| = \sum_{j=1}^K i_j$.
	Let $p_{\min} = \min_i p_i$. 
	Then we have
	$$
	1 - \mu(A_n) \le \sum_{j=1}^K (1-p_i)^n \le K(1-p_{\min})^n 
	\le c^{-1}\epsilon^{-s}\left(1-c\epsilon^s\right)^n.
	$$
	Note that if two points are in $X_i$,
	since $X_i \subset H_{\epsilon}(\textbf{z}_i)$,
	the distance between these two points is at most $\sqrt{d}\epsilon$.
	Thus, with probability at least $1 - c^{-1}\epsilon^{-s}\left(1-c\epsilon^s\right)^n$
	over iid $n$ samples $\Omega_n=\{\x_i\}_{i=1}^n$, 
	for any $\x \in X$, 
	there exists a point $\x_i$ such that $\|\x - \x_i\|_2 \le \sqrt{d}\epsilon$.
	By letting $\frac{1}{c^2\epsilon^{2s}} = n \ge \frac{1}{c\epsilon^s}$, 
	the proof is completed.
	Note that the probability in the above statement becomes   
	$1 - \sqrt{n}(1-\frac{1}{\sqrt{n}})^n$
	that goes to 1 as $n\to \infty$.
\end{proof}

\begin{proof}[Proof of Theorem~\ref{thm:gen}]
	Let $\mathcal{T}_r^{m_r}=\{\x_r^i\}_{i=1}^{m_r}$ be iid samples from $\mu_r$ on $U$
	and $\mathcal{T}_b^{m_b} = \{\x^i_b\}_{i=1}^{m_b}$ be 
	iid samples from $\mu_b$ on $\Gamma$.
	
	By Lemma~\ref{cor-sampling},
	with probability at least   
	\begin{equation} \label{app:thm-prob}
	(1 - \sqrt{m_r}(1-1/\sqrt{m_r})^{m_r})
	(1 - \sqrt{m_{b}}(1-1/\sqrt{m_{b}})^{m_{b}}),
	\end{equation}
	$\forall \x_r \in U$ and $\forall \x_b \in \Gamma$,
	there exists $\x_r' \in \mathcal{T}_r^{m_r}$
	and $\x_b' \in \mathcal{T}_b^{m_b}$ such that 
	$\|\x_r - \x_r'\| \le \sqrt{d}c_r^{-\frac{1}{d}}m_r^{-\frac{1}{2d}}$ and $\|\x_{b} - \x_{b}'\| \le \sqrt{d}c^{-\frac{1}{d-1}}_{b}m_{b}^{-\frac{1}{2(d-1)}}$.
	By letting $\epsilon_r = \sqrt{d}c_r^{-\frac{1}{d}}m_r^{-\frac{1}{2d}}$
	and $\epsilon_{b} = \sqrt{d}c^{-\frac{1}{d-1}}_{b}m_{b}^{-\frac{1}{2(d-1)}}$, 
	it follows from Lemma~\ref{app:lemma-lip}
	that with probability at least \eqref{app:thm-prob},
    \begin{equation*}
	\begin{split}
    \text{Loss}^{\text{PINN}}(h;\bm{\lambda})
	\le
	C_{\bm{m}} \cdot \left[\text{Loss}_{\bm{m}}^{\text{PINN}}(h;\bm{\lambda})
	+\lambda_{r,\bm{m}}^R \cdot \big[\mathcal{L}[h]\big]_{\alpha;U}^2
	+ \lambda_{b,\bm{m}}^R \cdot \big[\mathcal{B}[h]\big]_{\alpha;\Gamma}^2\right] + Q,
	\end{split}
	\end{equation*}
    where 
    $Q =3\lambda_r \sqrt{d}^{2\alpha} c_r^{-\frac{2\alpha}{d}} \big[f\big]_{\alpha;U}^2 m_r^{-\frac{\alpha}{d}} +3\lambda_{b} \sqrt{d}^{2\alpha} c_{b}^{-\frac{2\alpha}{d-1}}\big[g\big]_{\alpha;\Gamma}^2 m_{b}^{-\frac{\alpha}{d-1}}$,
    and 
    \begin{align*}
        C_{\bm{m}} &= 3\max\{\frac{C_r}{c_r}\sqrt{d}^{d} m_r^{\frac{1}{2}},  \frac{C_{b}}{c_{b}} \sqrt{d}^{d-1} m_{b}^{\frac{1}{2}}\},
        \quad 
        \lambda_{r,\bm{m}}^R = \frac{3\lambda_r \sqrt{d}^{2\alpha}c_r^{-\frac{2\alpha}{d}}m_r^{-\frac{\alpha}{d}}}{C_{\bm{m}}},
        \quad
        \lambda_{b,\bm{m}}^R = \frac{3\lambda_{b} \sqrt{d}^{2\alpha} c_{b}^{-\frac{2\alpha}{d-1}} m_{b}^{-\frac{\alpha}{d-1}}}{C_{\bm{m}}}.
    \end{align*}
    Let $C' = 3\max\{\lambda_r \sqrt{d}^{2\alpha} c_r^{-\frac{2\alpha}{d}}\big[f\big]_{\alpha;U}^2, \lambda_{b}\sqrt{d}^{2\alpha} c_{b}^{-\frac{2\alpha}{d-1}}\big[g\big]_{\alpha;\Gamma}^2  \}$.
    Then, we have
    \begin{align*}
        \text{Loss}^{\text{PINN}}(h;\bm{\lambda})
        \le C_{\bm{m}}\cdot \text{Loss}_{\bm{m}}(h;\bm{\lambda}, \bm{\lambda}_{\bm{m}}^R)
        +C'(m_r^{-\frac{\alpha}{d}} +  m_{b}^{-\frac{\alpha}{d-1}}),
    \end{align*}
    where $\bm{\lambda}_{\bm{m}}^R = (\lambda_{r,\bm{m}}^R,\lambda_{b,\bm{m}}^R)$.
    Then, the proof is completed.
\end{proof}

\section{Proof of Theorem~\ref{thm:conv-loss}} \label{app:thm:conv-loss}

\begin{proof}[Proof of Theorem~\ref{thm:conv-loss}]
	Suppose $m_r = \mathcal{O}(m_{b}^{\frac{d}{d-1}})$. 
	It then can be checked that 
	$\hat{\lambda}_{r,\bm{m}}^R = \hat{\lambda}_{b,\bm{m}}^R =\mathcal{O}(m_r^{-\frac{1}{2}-\frac{\alpha}{d}})$,
	where $\hat{\lambda}_{r,\bm{m}}$ and $\hat{\lambda}_{b,\bm{m}}$ are defined in \eqref{def:lambdas}.
	Let $\bm{\lambda}$ be a vector independent of $\bm{m}$
	and $\bm{\lambda}_{\bm{m}}^R$ be a vector satisfying \eqref{def:lambda-condition}.
	
	Let $h_{m} \in \mathcal{H}_{m}$ be a function that minimizes $\text{Loss}_{\bm{m}}(\cdot;\bm{\lambda}, \bm{\lambda}_{\bm{m}}^R)$.
	First, note that since $u_{\bm{m}}^* \in \mathcal{H}_{\bm{m}}$, 
	\begin{equation} \label{app:eqn-loss-convg}
	\begin{split}
	&
	\min\{\lambda_{r,\bm{m}}^R,\lambda_{b,\bm{m}}^R\} \left(R_r(h_{m}) + R_{b}(h_{m})\right) 
	\le \text{Loss}_{m}(h_{m};\bm{\lambda}, \bm{\lambda}_{\bm{m}}^R)
	\\
	&\le 
	\text{Loss}_{m}(u_{\bm{m}}^*;\bm{\lambda}, \bm{\lambda}_{\bm{m}}^R)
	\le
	\text{Loss}_{m}(u_{\bm{m}}^*;\bm{\lambda}, \bm{0}) + 
	\|\bm{\lambda}_{\bm{m}}^R\|_{\infty}\left(R_r(u_{\bm{m}}^*) + R_{b}(u_{\bm{m}}^*)\right). 
	\end{split}
	\end{equation}
	Since $\bm{\lambda}_{\bm{m}}^R \ge \hat{\bm{\lambda}}_{\bm{m}}^R$
	and $\|\bm{\lambda}_{\bm{m}}^R\|_\infty = \mathcal{O}(\|\hat{\bm{\lambda}}_{\bm{m}}^R\|_\infty)$, we have 
	$\frac{\max_j (\bm{\lambda}_{\bm{m}}^R)_j}{\min_j (\bm{\lambda}_{\bm{m}}^R)_j} =\mathcal{O}(1)$.
	Let $R^* = \sup_{\bm{m}} (R_r(u_{\bm{m}}^*) + R_b(u_{\bm{m}}^*))$.
	By the third assumption in \ref{assumption:convergence}, we have $R^* < \infty$.
	
	The second assumption of \ref{assumption:convergence}
    can be relaxed to the following condition:
    \begin{itemize}
        \item For each $m_r$, $\mathcal{H}_{m_r}$ contains a function $u_{m_r}^*$
		satisfying 
		$\text{Loss}_{m_r}^{\text{PINN}}(u_{m_r}^*;\bm{\lambda}) = \mathcal{O}(m_r^{-\frac{1}{2}-\frac{\alpha}{d}})$.
    \end{itemize}
	Let us write $m_r$ as $m$ for the sake of simplicity.
	We then have $R_r(h_{m}), R_{b}(h_{m})\le \mathcal{O}(R^*)$ for all $m$.
	Since $R_r(h_{m})=
	\big[\mathcal{L}[h_{m}]\big]_{\alpha;U}^2$
	and $R_{b}(h_{m})=\big[\mathcal{B}[h_{m}]\big]_{\alpha;U}^2$,
	the H\"{o}lder coefficients of $\mathcal{L}[h_{m}]$ and $\mathcal{B}[h_{m}]$ are uniformly bounded above.
	With the first assumption of \eqref{assumption:convergence},
    $\{\mathcal{L}[h_{m}]\}$ and $\{\mathcal{B}[h_{m}]\}$
    are uniformly bounded and uniformly equicontinuous sequences of functions in
    $C^{0,\alpha}(U)$ and $C^{0,\alpha}(\Gamma)$, respectively.
	By invoking the Arzela-Ascoli Theorem, 
	there exists a subsequence $h_{m_j}$
	and functions $G \in C^{0,\alpha}(U)$ and $B \in C^{0,\alpha}(\Gamma)$
	such that $\mathcal{L}[h_{m_j}] \to G$ 
	and $\mathcal{B}[h_{m_j}] \to B$ in $C^{0}(U)$ and $C^{0}(\Gamma)$, respectively, as $j \to \infty$.

    Since $\text{Loss}_m(h_m;\bm{\lambda},\bm{\lambda}_{\bm{m}}^R) = \mathcal{O}(m_r^{-\frac{1}{2}-\frac{\alpha}{d}})$,
	by combining it with Theorem~\ref{thm:gen}, we have that 
	with probability at least $(1 - \sqrt{m_r}(1-c_r/\sqrt{m_r})^{m_r})(1 - \sqrt{m_{b}}(1-c_{b}/\sqrt{m_{b}})^{m_{b}})$,
	\begin{equation*}
	\text{Loss}(h_{m};\bm{\lambda},\bm{0}) = \mathcal{O}(m^{-\frac{\alpha}{d}}). 
	\end{equation*} 
	Hence, the probability of
	$\lim_{m \to \infty} 
	\text{Loss}(h_{m};\bm{\lambda},\bm{0}) = 0$
	is one. Thus, with probability 1, 
	\begin{align*}
	0 &= \lim_{j \to \infty} \text{Loss}(h_{m_j};\bm{\lambda},\bm{0})
	\\
	&= \lim_{j \to \infty} \lambda_r \int_{U} \|\mathcal{L}[h_{m_j}](\x_r) - f(\x_r) \|^2d\mu_r(\x_r) +  \lambda_{b}\int_{\Gamma} \|\mathcal{B}[h_{m_j}](\x_{b}) - g(\x_{b}) \|^2 d\mu_{b}(\x_{b})
	\\
	&=\lambda_{f}\int_{U} \|G(\x_r) - f(\x_r) \|^2d\mu_r(\x_r) + 
	\lambda_{b} \int_{\Gamma} \|B_k(\x_{b}) - g(\x_{b}) \|^2 d\mu_{b}(\x_{b}),
	\end{align*}
	which shows that 
	$G = f$ in $L^2(U;\mu_r)$ and $B = g$ in $L^2(\Gamma;\mu_{b})$.
	Note that since $\mathcal{L}[h_{m_j}]$ and $\mathcal{B}[h_{m_j}]$ are uniformly bounded above and uniformly converge to $G$ and $B$, respectively,
	the third equality holds by Lebesgue's Dominated Convergence Theorem.
	Since the subsequence was arbitrary,
	we conclude that $\mathcal{L}[h_{m}] \to f$ in $L^2(U;\mu_r)$ and $\mathcal{B}[h_m] \to g$
	in $L^2(\Gamma;\mu_{b})$ as $m \to \infty$.
	Furthermore, since $\{\mathcal{L}[h_{m}]\}$ and $\{\mathcal{B}[h_m]\}$ are equicontinuous,
	its convergence mode is improved to $C^{0}$.
\end{proof}

\section{Proof of Theorem~\ref{thm:main-elliptic}} \label{app:thm:main-elliptic}
For readability, 
the existence and the uniqueness of the classical solution to \eqref{def:elliptic-PDE} is stated as follow.
\begin{theorem}[Theorem 6.13 of \cite{Gilbarg_15_EllipticPDEs}] \label{thm:elliptic-existence}
	For $0 < \alpha < 1$,
	let $U$ satisfy an exterior sphere condition at every boundary point.
	Let the operator $\mathcal{L}$ be strictly elliptic in $U$
	with coefficients in $C^\alpha(U)$ and $\tilde{c}(\x) \le 0$.
	Then, the Dirichlet problem of \eqref{def:elliptic-PDE}
	has a unique solution in $C^{2,\alpha}(U)\cap C^{0}(\bar{U})$
	for all $f \in C^\alpha(U)$
	and all $g \in C^{0}(\partial U)$.
\end{theorem}

\begin{lemma} \label{lemma:elliptic-stability}
	Suppose Assumption~\ref{assumption:elliptic} holds and $\tilde{c}(\x) \le 0$.
	Then, 
	there exists the unique classical solution $u^*$ to the PDE \eqref{def:elliptic-PDE}.
	Furthermore, there exists a positive constant $C$ that depends only on $U$, $\tilde{\bm{b}}(\x)$, $\lambda_0$ and $v$
	such that 
	\begin{equation}
	\|u^* - h\|_{C^0(U)} \le C(\|f - \mathcal{L}[h]\|_{C^0(U)} + \|g - h\|_{C^0(\partial U)}),
	\end{equation}
	for any $h \in C^{0}(U)$
	satisfying $\mathcal{L}[h] \in C^{0,\alpha}(U)$
	and $h \in C^{0,\alpha}(\partial U)$.
\end{lemma}
\begin{proof}[Proof of Lemma~\ref{lemma:elliptic-stability}]
	By Theorem~\ref{thm:elliptic-existence},
	there exists a unique classical solution $\xi \in C^{2,\alpha}(U)\cap C^{0}(\bar{U})$
	to
	\begin{equation*} 
	\begin{split}
	\mathcal{L}[\xi] = 0, \text{ on } U, \qquad
	\xi|_{\partial U} = h|_{\partial U} - g.
	\end{split}
	\end{equation*}
	Since $\xi \in C^{2,\alpha}(U)\cap C^{0}(\bar{U})$, $\tilde{c}(\x) \le 0$ 
	and Assumption~\ref{assumption:elliptic} holds, 
	one can apply the Weak Maximum principle (e.g. Corollary 3.11 in \cite{Gilbarg_15_EllipticPDEs}) to have  
	$$
	\|\xi\|_{C^{0}(U)} \le \|h - g\|_{C^0(\partial U)}.
	$$

	By Theorem~\ref{thm:elliptic-existence}, the PDE \eqref{def:elliptic-PDE} has a unique classical solution $u^*$.
	Let $\tilde{e} =u^* - h + \xi$.
	Then, $\tilde{e}$ satisfies 
	$\mathcal{L}[\tilde{e}] = f - \mathcal{L}[h]$ in $U$
	and $\tilde{e}|_{\partial U} = 0$.
	From an apriori bound of the elliptic equation (e.g. Theorem 3.7. of \cite{Gilbarg_15_EllipticPDEs}), we have
	$$
	\|\tilde{e}\|_{C^0(U)} \le \frac{C}{\lambda_0}\|f - \mathcal{L}[h]\|_{C^0(U)},
	$$
	where $C$ is a constant depending only on $\text{diam}(U)$ and $\sup_{U} \frac{\|\tilde{\bm{b}}(\x)\|}{\lambda_0}$.
	Hence,
	\begin{align*}
	\|u^* - h\|_{C^0(U)} \le \|\tilde{e}\|_{C^{0}(U)} + \|\xi\|_{C^0(U)}
	\le \frac{C}{\lambda_0}\|f - \mathcal{L}[h]\|_{C^0(U)}
	+ \|h - g\|_{C^0(\partial U)},
	\end{align*}
	which completes the proof.
\end{proof}

\begin{proof}[Proof of Theorem~\ref{thm:main-elliptic}]
	Let $\{h_m\}_{m\ge 1}$ be the sequence of the minimizers from Theorem~\ref{thm:conv-loss}.
	It follows from Theorem~\ref{thm:conv-loss}
	that $\|h_m - g\|_{C^0(\partial U)} \to 0$ 
	and $\|f - \mathcal{L}[h_m]\|_{C^0(U)} \to 0$.
	By Assumption~\ref{assumption:convergence}, 
	$\mathcal{L}[h] \in C^{0,\alpha}(U)$
	and $h \in C^{0,\alpha}(\partial U)$ for all $h \in \mathcal{H}_{\bm{m}}$.
	Lemma~\ref{lemma:elliptic-stability} shows 
	$h_m \to u^*$ in $C^0(U)$,
	which completes the proof.
\end{proof}

\section{Proof of Theorem~\ref{thm:elliptic-H1}} \label{app:thm:elliptic-H1}
\begin{proof}
	The proof continues from the proof of Lemma~\ref{lemma:elliptic-stability}. 
	Since $\tilde{e}|_{\partial U} = 0$, 
	it follows from Assumption~\ref{assumption:elliptic} that 
	\begin{align*}
	-\langle \mathcal{L}[\tilde{e}], \tilde{e} \rangle_{L^2(U)}
	\ge \frac{\lambda_0}{2} \|D\tilde{e}\|_{L^2(U)}^2 - \lambda_0 v^2 \|\tilde{e}\|_{L^2(U)}^2.
	\end{align*}
	Thus, 
	\begin{align*}
	\|D\tilde{e}\|_{L^2(U)}^2 &\le 
	\|\mathcal{L}[\tilde{e}]\|_{L^2(U)}
	\|\tilde{e}\|_{L^2(U)} + 2v^2 \|\tilde{e}\|_{L^2(U)}^2 \\
	&\le C'\left[\|\mathcal{L}[\tilde{e}]\|^2_{L^2(U)} + \|\tilde{e}\|_{L^2(U)}^2\right],
	\end{align*}
	where $C'$ is a constant independent of $h$.
	Therefore, 
	$$
	\|\tilde{e}\|_{H^1(U)}^2 
	= \|D\tilde{e}\|_{L^2(U)}^2 + \|\tilde{e}\|_{L^2(U)}^2 
	\le C^{''}\left[\|\mathcal{L}[\tilde{e}]\|^2_{L^2(U)} + \|\tilde{e}\|_{L^2(U)}^2\right],
	$$
	for some constant $C^{''}$.
	Then,
	\begin{align*}
	\|u^* - h\|^2_{H^1(U)} &\le 2\|\tilde{e}\|^2_{H^1(U)} + 2\|\xi\|^2_{H^1(U)} \\
	&\le C^{'''} \left(\|f - \mathcal{L}[h]\|_{C^0(U)} + \|g - h\|_{C^0(\partial U)}\right)^2 + 2\|\xi\|^2_{H^1(U)},
	\end{align*}
	for some constant $C^{'''}$.
	Since $\xi = 0$ and $h_m = g$ on $\partial U$,
	the proof is completed by Theorem~\ref{thm:conv-loss} 
	as $\|f - \mathcal{L}[h_m]\|_{C^0(U)} \to 0$.
\end{proof}


\section{Proof of Theorem~\ref{thm:main-parabolic}} \label{app:thm:main-parabolic}
For reader's convenience, let us recall a result on the existence of the classical solution to the parabolic equation.
\begin{theorem}[Corollary 2 in Chapter 3 of \cite{Friedman_08_ParabolicPDEs}] \label{thm:parabolic-existence}
	Suppose Assumption~\eqref{assumption:parabolic} holds.
	suppose further that for every point $x' \in \partial U$, there exists a closed ball $\bar{B}$ in $\mathbb{R}^n$ such that 
	$\bar{B}\cap \overline{U} = \{x'\}$.
	Then, for any continuous function $\varphi$ on $\partial U \times (0,T]$ and $g$ on $
	\overline{U}\times \{t=0\}$
	such that $g(\x) = \varphi(\x, 0)$ on $\partial U$,
	there exists a unique classical solution $u^*$ to \eqref{def:Parabolic}.
\end{theorem}

\begin{lemma} \label{lem:stability-parabolic}
	Suppose Assumption~\ref{assumption:convergence}
	and ~\ref{assumption:parabolic} hold
	and $\tilde{c}(\x,t) \le 0$.
	Let $u^*$ be the solution to the PDE \eqref{def:Parabolic}
	and let $\psi(\x,t)$ be the combined boundary and initial conditions on 
	$\Gamma_T = \overline{U}_T - U_T$.
	Then, there exists a positive universal constant $C$ that depends only on
	the domain and the PDE \eqref{def:Parabolic} such that 
	\begin{equation*}
	\|h - u^*\|_{C^0(\Omega)} \le C\left[ \|f + h_t -\mathcal{L}[h]\|_{C^0(\Omega)}
	+ \|h - \psi\|_{C^0(\Gamma_T)}
	\right],
	\end{equation*}
	for any $h \in C_{2+\alpha} \cap C^{0}(\overline{U}\times [0,T])$.
\end{lemma}
\begin{proof}
	Let us combine the boundary and the initial conditions into one condition:
	$$
	u(x,t) = \psi(x,t) \quad \text{on } \Gamma_T = \overline{U}_T - U_T.
	$$
	For a fixed $h \in C_{2+\alpha} \cap C^{0}(\overline{U}\times [0,T])$, let $\xi$ be the solution to
	\begin{equation*} 
	\begin{cases}
	-\xi_t + \mathcal{L}[\xi] = 0,  & \text{in } U_T \\
	\xi = \psi - h, & \text{on } \Gamma_T
	\end{cases}
	\end{equation*}
	We note that since $\psi - h$ is continuous on $\Gamma_T$, the existence of $\xi$ is guaranteed by Theorem~\ref{thm:parabolic-existence}.
	Note that since $\xi$ is continuous, and $\mu(\Omega) > 0$,
	it follows from the Weak Maximum Principle that
	$$
	\|\xi\|_{C^0(\Omega)} \le \|\psi - h\|_{C^0(\Gamma_T)}.
	$$
	Let $u^*$ be a unique classical solution to the PDE \eqref{def:Parabolic}, 
	and let $\tilde{e} = u^* - h - \xi$.
	Since Assumption~\ref{assumption:parabolic} holds
	and $\tilde{c}(\x,t) \le 0$,
	it follows from a consequence of the weak maximum principle (e.g. p 42, Chapter 2 of Friedman \cite{Friedman_08_ParabolicPDEs}) that 
	\begin{align*}
	\|\tilde{e}\|_{C^0(\Omega)} \le C' \|f + h_t - \mathcal{L}[h]\|_{C^0(\Omega)},
	\end{align*}
	for some constant $C'$.
	Therefore,
	\begin{align*}
	\|h - u^*\|_{C^0(\Omega)} &\le 
	\|\tilde{e}\|_{C^0(\Omega)} + \|\xi\|_{C^0(\Omega)}
	\\
	&\le 
	C''\left[\|f + h_t - \mathcal{L}[h]\|_{C^0(\Omega)}
	+
	\|\psi - h\|_{C^0(\Gamma_T)}
	\right].
	\end{align*}
	Here $C''$ is a positive universal constant that is independent of $h$.
\end{proof}	

\begin{proof}[Proof of Theorem~\ref{thm:main-parabolic}]
	Let $\{h_m\}_{m\ge 1}$ be the sequence of minimizers of \eqref{def:problem}.
	Note that from Theorem~\ref{thm:conv-loss}, we have
	\begin{align*}
	\lim_{m\to \infty} \|(h_m)_t - \mathcal{L}[h_m] + f\|_{C^0(\Omega)} = 0,
	\qquad
	\lim_{m\to \infty} \|h_m - \psi\|_{C^0(\Gamma_T)} = 0.
	\end{align*}
	By Lemma~\ref{lem:stability-parabolic},
	it readily follows that 
	$\lim_{m \to \infty} \|h_m - u^*\|_{C^0(\Omega)}  =0$.
\end{proof}	
	
\section{Proof of Theorem~\ref{thm:parabolic-H1}} \label{app:thm:parabolic-H1}
\begin{proof}	
	The proof continues from the proof of Lemma~\ref{lem:stability-parabolic}.
	Note that 
	$\frac{d}{dt}\tilde{e}_m - \mathcal{L}[\tilde{e}_m] = f + (h_m)_t - \mathcal{L}[h_m]$ on $U_T$ and 
	$\tilde{e}_m|_{\Gamma_T} = 0$.
	For a fixed $t$, by the integration by parts, one obtains 
	\begin{align*}
	\langle (\tilde{e}_m)_t - \mathcal{L}[\tilde{e}_m], \tilde{e}_m \rangle_{L^2(U)} 
	&= \frac{1}{2}\frac{d}{dt} \|\tilde{e}_m\|_{L^2(U)}^2
	- \langle \mathcal{L}[\tilde{e}_m], \tilde{e}_m \rangle_{L^2(U)} 
	\\
	&\ge \frac{1}{2}\frac{d}{dt} \|\tilde{e}_m(t)\|_{L^2(U)}^2
	+ \beta \|\tilde{e}_m(t)\|_{H^1(U)}^2 - \gamma \|\tilde{e}_m(t)\|_{L^2(U)}^2,
	\end{align*}
	for some $\beta > 0$ and $\gamma \ge 0$.
	Also, note that 
	$$
	\langle (\tilde{e}_m)_t - \mathcal{L}[\tilde{e}_m], \tilde{e}_m \rangle_{L^2(U)} 
	\le \frac{1}{2}\| f +(h_m)_t - \mathcal{L}[h_m]\|_{L^2(U)}^2 + \frac{1}{2}\|\tilde{e}_m(t)\|_{L^2(U)}^2.
	$$
	By combining the above two inequalities, we have
	\begin{align*}
	\frac{d}{dt} \|\tilde{e}_m(t)\|_{L^2(U)}^2 + 
	\beta' \|\tilde{e}_m(t)\|_{H^1(U)}^2
	\le C_1\|\tilde{e}_m(t)\|_{L^2(U)}^2 + C_2\| f +(h_m)_t - \mathcal{L}[h_m]\|_{L^2(U)}^2,
	\end{align*}
	for some positive constants $\beta', C_1, C_2$ that are independent of $m$.
	By integrating it with respect to $t$ from $0$ to $T$,
	\begin{align*}
	&\|\tilde{e}_m(T)\|_{L^2(U)}^2 - \|\tilde{e}_m(0)\|_{L^2(U)}^2
	+ \beta' \|\tilde{e}_m\|_{L^2(0,T;H^1(U))}^2 \\
	&\le C_1\|\tilde{e}_m\|_{L^2(0,T;L^2(U))}^2 + C_2\| f +(h_m)_t - \mathcal{L}[h_m]\|_{L^2(0,T;L^2(U))}^2.
	\end{align*}
	Since $\tilde{e}_m|_{\Gamma_T} = 0$, 
	we have $\|\tilde{e}_m(0)\|_{L^2(U)}^2=0$.
	Hence, 
	\begin{align*}
	\|\tilde{e}_m\|_{L^2(0,T;H^1(U))}^2
	\le C'_1\|\tilde{e}_m\|_{L^2(0,T;L^2(U))}^2 + C'_2\| f +(h_m)_t - \mathcal{L}[h_m]\|_{L^2(0,T;L^2(U))}^2.
	\end{align*}

	By Lemma~\ref{lem:stability-parabolic}, we have
	\begin{align*}
	\|\tilde{e}_m\|_{L^2(0,T;H^1(U))}^2
	\le C''\|f +(h_m)_t - \mathcal{L}[h_m]\|_{C^0(0,T;C^0(U))}^2,
	\end{align*}
	for some constant $C''$.
	Hence, 
	$\lim_{m\to\infty} \|\tilde{e}_m\|_{L^2(0,T;H^1_0(U))} = 0$.
	Therefore, $h_m \to u^*$ in $L^2(0,T;H^1(U))$.
\end{proof}

%% file: acknowledgement.tex
\section*{Acknowledgement}
The authors would like to thank Dr. Hongjie Dong,
Dr. Seick Kim
for the helpful discussion on the Schauder approch.
G. E. Karniadakis acknowledges support by the PhILMS grant DE-SC0019453.
Furthermore, the authors would like to thank anonymous referees
for their invaluable comments improving an early version of this manuscript.

%% file: main.bbl
\begin{thebibliography}{10}

\bibitem{Attouch_14Book_variational}
H.~Attouch, G.~Buttazzo, and G.~Michaille.
\newblock {\em Variational analysis in Sobolev and BV spaces: applications to
  PDEs and optimization}.
\newblock SIAM, 2014.

\bibitem{Baker_19_DCworkshop}
N.~Baker, F.~Alexander, T.~Bremer, A.~Hagberg, Y.~Kevrekidis, H.~Najm,
  M.~Parashar, A.~Patra, J.~Sethian, S.~Wild, and ET~AL.
\newblock Workshop report on basic research needs for scientific machine
  learning: Core technologies for artificial intelligence.
\newblock Technical report, USDOE Office of Science (SC), Washington, DC
  (United States), 2019.

\bibitem{Baydin_17AD}
A.~G. Baydin, B.~A. Pearlmutter, A.~A. Radul, and J.~M. Siskind.
\newblock Automatic differentiation in machine learning: a survey.
\newblock {\em The Journal of Machine Learning Research}, 18(1):5595--5637,
  2017.

\bibitem{Berg_18_Unified}
J.~Berg and K.~Nystr{\"o}m.
\newblock A unified deep artificial neural network approach to partial
  differential equations in complex geometries.
\newblock {\em Neurocomputing}, 317:28--41, 2018.

\bibitem{bottou2008tradeoffs}
L.~Bottou and O.~Bousquet.
\newblock The tradeoffs of large scale learning.
\newblock In {\em Advances in Neural Information Processing Systems}, pages
  161--168, 2008.

\bibitem{Brezis_10Book_functional}
H.~Brezis.
\newblock {\em Functional analysis, Sobolev spaces and partial differential
  equations}.
\newblock Springer Science \& Business Media, 2010.

\bibitem{Calder_19_consistency}
J.~Calder.
\newblock Consistency of {L}ipschitz learning with infinite unlabeled data and
  finite labeled data.
\newblock {\em SIAM Journal on Mathematics of Data Science}, 1(4):780--812,
  2019.

\bibitem{Darbon_19_NNsHJ}
Jerome Darbon, Gabriel~P. Langlois, and Tingwei Meng.
\newblock Overcoming the curse of dimensionality for some hamilton--jacobi
  partial differential equations via neural network architectures.
\newblock {\em arXiv preprint arXiv:1910.09045}, 2019.

\bibitem{Darbon_20_NNsHJ}
Jerome Darbon and Tingwei Meng.
\newblock On some neural network architectures that can represent viscosity
  solutions of certain high dimensional {H}amilton--{J}acobi partial
  differential equations.
\newblock {\em arXiv preprint arXiv:2002.09750}, 2020.

\bibitem{Dissanayake_94_ANN-PDE}
M.~Dissanayake and N.~Phan-Thien.
\newblock Neural-network-based approximations for solving partial differential
  equations.
\newblock {\em Communications in Numerical Methods in Engineering},
  10(3):195--201, 1994.

\bibitem{Weinan_18DeepRitz}
W.~E and B.~Yu.
\newblock The deep {R}itz method: a deep learning-based numerical algorithm for
  solving variational problems.
\newblock {\em Communications in Mathematics and Statistics}, 6(1):1--12, 2018.

\bibitem{Evans_15_measure}
L.~C. Evans and R.~F. Gariepy.
\newblock {\em Measure theory and fine properties of functions}.
\newblock CRC press, 2015.

\bibitem{Finlay_18_Lipschitz}
C.~Finlay, J.~Calder, B.~Abbasi, and A.~Oberman.
\newblock Lipschitz regularized deep neural networks generalize and are
  adversarially robust.
\newblock {\em arXiv preprint arXiv:1808.09540}, 2018.

\bibitem{Friedman_08_ParabolicPDEs}
A.~Friedman.
\newblock {\em Partial differential equations of parabolic type}.
\newblock Courier Dover Publications, 2008.

\bibitem{Gilbarg_15_EllipticPDEs}
D.~Gilbarg and N.~S. Trudinger.
\newblock {\em Elliptic partial differential equations of second order}.
\newblock Springer, 2015.

\bibitem{Glorot2010understanding}
X.~Glorot and Y.~Bengio.
\newblock Understanding the difficulty of training deep feedforward neural
  networks.
\newblock In {\em Proceedings of the thirteenth international conference on
  artificial intelligence and statistics}, pages 249--256, 2010.

\bibitem{Grohs_18_BSPDEs}
P.~Grohs, F.~Hornung, A.~Jentzen, and P.~Von Wurstemberger.
\newblock A proof that artificial neural networks overcome the curse of
  dimensionality in the numerical approximation of black-scholes partial
  differential equations.
\newblock {\em arXiv preprint arXiv:1809.02362}, 2018.

\bibitem{HanE_18_DLSPDEs}
J.~Han, A.~Jentzen, and W.~E.
\newblock Solving high-dimensional partial differential equations using deep
  learning.
\newblock {\em Proceedings of the National Academy of Sciences},
  115(34):8505--8510, 2018.

\bibitem{Houska_19_Global}
B.~Houska and B.~Chachuat.
\newblock Global optimization in hilbert space.
\newblock {\em Mathematical programming}, 173(1-2):221--249, 2019.

\bibitem{Jagtap_19LAAF}
A.~D. Jagtap, K.~Kawaguchi, and G.~E. Karniadakis.
\newblock Locally adaptive activation functions with slope recovery term for
  deep and physics-informed neural networks.
\newblock {\em arXiv preprint arXiv:1909.12228}, 2019.

\bibitem{Kingma_14_Adam}
D.~P. Kingma and J.~Ba.
\newblock Adam: A method for stochastic optimization.
\newblock {\em arXiv preprint arXiv:1412.6980}, 2014.

\bibitem{Lagaris_98_ANN-ODE-PDE}
I.~E. Lagaris, A.~Likas, and D.~I. Fotiadis.
\newblock Artificial neural networks for solving ordinary and partial
  differential equations.
\newblock {\em IEEE transactions on Neural Networks}, 9(5):987--1000, 1998.

\bibitem{Lagaris_00_ANN-Irregular}
I.~E. Lagaris, A.~C. Likas, and G.~D. Papageorgiou.
\newblock Neural-network methods for boundary value problems with irregular
  boundaries.
\newblock {\em IEEE Transactions on Neural Networks}, 11(5):1041--1049, 2000.

\bibitem{Lecun_Nature15_DeepLearning}
Y.~LeCun, Y.~Bengio, and G.~Hinton.
\newblock Deep learning.
\newblock {\em Nature}, 521(7553):436--444, 2015.

\bibitem{Liu_89_LBFGS}
D.~C. Liu and J.~Nocedal.
\newblock On the limited memory bfgs method for large scale optimization.
\newblock {\em Mathematical programming}, 45(1-3):503--528, 1989.

\bibitem{Lu_19_Deepxde}
L.~Lu, X.~Meng, Z.~Mao, and G.~E. Karniadakis.
\newblock Deepxde: A deep learning library for solving differential equations.
\newblock {\em arXiv preprint arXiv:1907.04502}, 2019.

\bibitem{Mao_20_HighSpeedFlows}
Z.~Mao, A.~D. Jagtap, and G.~E. Karniadakis.
\newblock Physics-informed neural networks for high-speed flows.
\newblock {\em Computer Methods in Applied Mechanics and Engineering},
  360:112789, 2020.

\bibitem{Mohri_18_FoundationsML}
M.~Mohri, A.~Rostamizadeh, and A.~Talwalkar.
\newblock {\em Foundations of machine learning}.
\newblock MIT press, 2018.

\bibitem{niyogi1999generalization}
P.~Niyogi and F.~Girosi.
\newblock Generalization bounds for function approximation from scattered noisy
  data.
\newblock {\em Advances in Computational Mathematics}, 10(1):51--80, 1999.

\bibitem{Pang_SISC19_fPINNs}
G.~Pang, L.~Lu, and G.~E. Karniadakis.
\newblock f{PINN}s: Fractional physics-informed neural networks.
\newblock {\em SIAM Journal on Scientific Computing}, 41(4):A2603--A2626, 2019.

\bibitem{Pinkus_99_ATofMLP}
A.~Pinkus.
\newblock Approximation theory of the mlp model in neural networks.
\newblock {\em Acta numerica}, 8:143--195, 1999.

\bibitem{Raissi_19_PINNs}
M.~Raissi, P.~Perdikaris, and G.~E. Karniadakis.
\newblock Physics-informed neural networks: A deep learning framework for
  solving forward and inverse problems involving nonlinear partial differential
  equations.
\newblock {\em Journal of Computational Physics}, 378:686--707, 2019.

\bibitem{Raissi_Nature20_HFM}
M.~Raissi, A.~Yazdani, and G.~E. Karniadakis.
\newblock Hidden fluid mechanics: Learning velocity and pressure fields from
  flow visualizations.
\newblock {\em Science}, 367(6481):1026--1030, 2020.

\bibitem{Ruder_16_GDoverview}
S.~Ruder.
\newblock An overview of gradient descent optimization algorithms.
\newblock {\em arXiv preprint arXiv:1609.04747}, 2016.

\bibitem{Sirignano_JCP18_DGM}
J.~Sirignano and K.~Spiliopoulos.
\newblock D{GM}: A deep learning algorithm for solving partial differential
  equations.
\newblock {\em Journal of Computational Physics}, 375:1339--1364, 2018.

\bibitem{Song_19_fPINNs}
F.~Song, G.~Pange, C.~Meneveau, and G.~E. Karniadakis.
\newblock Fractional physical-inform neural networks (f{PINN}s) for turbulent
  flows.
\newblock {\em Bulletin of the American Physical Society}, 2019.

\bibitem{Wang_20_GradPathPINN}
S.~Wang, Y.~Teng, and P.~Perdikaris.
\newblock Understanding and mitigating gradient pathologies in physics-informed
  neural networks.
\newblock {\em arXiv preprint arXiv:2001.04536}, 2020.

\bibitem{yosida1988functional}
K.~Yosida.
\newblock {\em Functional analysis}, volume 123.
\newblock springer, 1988.

\bibitem{Zhang_19_SPDE_PINNs}
D.~Zhang, L.~Guo, and G.~E. Karniadakis.
\newblock Learning in modal space: Solving time-dependent stochastic pdes using
  physics-informed neural networks.
\newblock {\em SIAM Journal on Scientific Computing}, 42(2):A639--A665, 2020.

\end{thebibliography}
